\documentclass[sn-mathphys,Numbered]{sn-jnl}% Math and Physical Sciences Reference Style
\usepackage{graphicx}%
\usepackage{multirow}%
\usepackage{amsmath,amssymb,amsfonts}%
\usepackage{amsthm}%
\usepackage{mathrsfs}%
\usepackage[title]{appendix}%
\usepackage{xcolor}%
\usepackage{tikz-cd}
\usepackage{textcomp}%
\usepackage{manyfoot}%
\usepackage{booktabs}%
\usepackage{algorithm}%
\usepackage{algorithmicx}%
\usepackage{algpseudocode}%
\usepackage{listings}%
\usepackage{comment}
\newcommand{\Spec}{\mathrm{Spec}}

\newcommand{\Cc}{\mathbb{C}}

\newcommand{\Qq}{\mathbb{Q}}

\newcommand{\Rr}{\mathbb{R}}

\newcommand{\Zz}{\mathbb{Z}}
\newcommand{\ZZ}{\mathbb{Z}}

\newcommand{\Center}{\operatorname{center}}

\newcommand{\Exc}{\operatorname{Exc}}

\newcommand{\rk}{\operatorname{rank}}

\newcommand{\mld}{{\operatorname{mld}}}
\newcommand{\num}{{\operatorname{num}}}
\newcommand{\tang}{{\operatorname{tang}}}
\newcommand{\pld}{{\operatorname{pld}}}

\newcommand{\Supp}{\operatorname{Supp}}

\newcommand{\mult}{\operatorname{mult}}

\newcommand{\Div}{\operatorname{Div}}

\newcommand{\lf}{\lfloor}
\newcommand{\rf}{\rfloor}

\newcommand{\Ff}{\mathcal{F}}

\newcommand{\Oo}{\mathcal{O}}
\newcommand{\Ii}{\Gamma}

\newcommand{\Sing}{\mathrm{Sing}}
\newcommand{\Pic}{\mathrm{Pic}}

\newtheorem{thm}{Theorem}[section]
\newtheorem{conj}[thm]{Conjecture}
\newtheorem{cor}[thm]{Corollary}
\newtheorem{lem}[thm]{Lemma}
\newtheorem{prop}[thm]{Proposition}

\theoremstyle{definition}
\newtheorem{defn}[thm]{Definition}

\theoremstyle{definition}
\newtheorem{rem}[thm]{Remark}

\newtheorem{ex}[thm]{Example}

\newtheorem{exrem}[thm]{Example-Remark}

\begin{document}

\title[Complements, index theorem, and minimal log discrepancies]{Complements, index theorem, and minimal log discrepancies of foliated surface singularities}


\author*[1]{\fnm{Jihao} \sur{Liu}}\email{jliu@northwestern.edu}

\author[2]{\fnm{Fanjun} \sur{Meng}}\email{fmeng3@jhu.edu}
\equalcont{These authors contributed equally to this work.}

\author[3]{\fnm{Lingyao} \sur{Xie}}\email{lingyao@math.utah.edu}
\equalcont{These authors contributed equally to this work.}

\affil*[1]{\orgdiv{Department of Mathematics}, \orgname{Northwestern University}, \orgaddress{\street{2033 Sheridan Road, \city{Evanston}, \postcode{60201}, \state{Illinois}, \country{USA}}}

\affil[2]{\orgdiv{Department of Mathematics,}, \orgname{Johns Hopkins University}, \orgaddress{\street{3400 N. Charles Street}, \city{Baltimore}, \postcode{21218}, \state{Maryland}, \country{USA}}}

\affil[3]{\orgdiv{Department of Mathematics,}, \orgname{The University of Utah}, \orgaddress{\street150 South 1400 East}, \city{Salt Lake City}, \postcode{84112}, \state{Utah}, \country{USA}}}


\abstract{We present an extension of several results on pairs and varieties to foliated surface pairs. We prove the boundedness of local complements, the local index theorem, and the uniform boundedness of minimal log discrepancies (mlds), as well as establishing the existence of uniform rational lc polytopes. Furthermore, we address two questions posed by P. Cascini and C. Spicer on foliations, providing negative responses. We also demonstrate that the Grauert-Riemenschneider type vanishing theorem generally fails for lc foliations on surfaces. In addition, we determine the set of minimal log discrepancies for foliated surface pairs with specific coefficients, which leads to the recovery of Y.-A. Chen's proof on the ascending chain condition conjecture for mlds for foliated surfaces.}




\maketitle

\section{Introduction}\label{sec1}

We work over the field of complex numbers $\mathbb C$.

The study of foliations is a major topic in birational geometry. In recent years, there has been significant progress on the minimal model program for foliated varieties in dimension $\leq 3$, as seen in \cite{CS20,Spi20,CS21,SS22}. While the global structures for foliations of dimension $\leq 3$ are mostly settled from the point of view of the minimal model program, there is still much to be explored regarding the local structure of foliations, i.e. the singularities of foliations.

A recent important contribution to the study of singularities of foliations is Y.-A. Chen's classification of $\Qq$-Gorenstein lc foliated surface singularities \cite[Theorem 0.1]{Che23}. As an application, Y.-A. Chen shows that the minimal log discrepancies (mlds) of foliated surface pairs with DCC coefficients satisfy the ascending chain condition (ACC) \cite[Theorem 0.2]{Che23}. It then becomes interesting to ask whether other standard conjectures in birational geometry, such as the boundedness of complements and Shokurov's local index conjecture, will also hold for foliations. In this paper, we study the analogues of these conjectures for surfaces.

\smallskip

\noindent\textbf{Local complements for foliated surfaces}. Complements theory is an essential tool in modern birational geometry. Birkar famously proved the existence of complements for any Fano type variety \cite[Theorem 1.8]{Bir19}, which was later used in the proof of the BAB conjecture \cite{Bir21}. This theory can be naturally applied to foliations as well. For example, foliated $1$-complements have played a crucial role in proving the existence of flips for rank $1$ foliated threefolds \cite{CS20}.

It is natural to ask whether the boundedness of complements also holds for foliations. In this paper, we prove the boundedness of local complements for foliated surfaces:

\begin{thm}\label{thm: lc complement foliated surface local}
Let $\epsilon$ be a non-negative real number and $\Ii\subset[0,1]$ be  DCC set. Then there exists a positive integer $N$ depending only on $\epsilon$ and $\Ii$ satisfying the following.

Assume that $(X\ni x,\Ff,B)$ is an $\epsilon$-lc foliated surface germ such that $B\in\Ii$. Then $(X\ni x,\Ff,B)$ has an $(\epsilon,N)$-complement $(X\ni x,\Ff,B^+)$, i.e. an $\epsilon$-lc foliated germ $(X\ni x,\Ff,B^+)$ such that $NB^+\geq N\lfloor B\rfloor+\lfloor(N+1)\{B\}\rfloor$ and $N(K_{\Ff}+B^+)$ is Cartier near $x$. 

Moreover, if $\bar\Ii\subset\mathbb Q$, then we may take $B^+\geq B$.
\end{thm}

Theorem \ref{thm: lc complement foliated surface local} provides positive evidence supporting the boundedness of complements for foliations. We refer the reader to Definition \ref{defn: complement} for a formal definition of complements for surfaces and to Conjecture \ref{conj: complement conjecture foliation} for a formal statement on the conjecture of the boundedness of complements for foliations. We plan to prove the global and relative cases of the boundedness of complements for foliated surfaces in future work.

In the special case when $\Ii={0}$, \cite[Question 4]{CS} inquires about whether a $1$-complement for a (relatively) Fano-type foliation always exists. For the local case of foliated surfaces, we have the following result, which provides a negative answer to  \cite[Question 4]{CS}:

\begin{thm}\label{thm: 1 2 complement foliation}
Let $(X\ni x,\Ff)$ be a foliated lc surface germ such that $\rk\Ff=1$. Then $(X\ni x,\Ff)$ has a $2$-complement. Moreover, there are cases when $(X\ni x,\Ff)$ do not have a $1$-complement.
\end{thm}
Theorem \ref{thm: 1 2 complement foliation} suggests that even for rank $1$ foliations on surfaces, the existence of a $1$-complement may be too optimistic, despite the expectation of boundedness of complements for foliations. We also note that non-exceptional surface singularities always have either a $1$-complement or a $2$-complement. Therefore, it is possible that the explicit values of $n$ in the boundedness of foliated $n$-complements are related to the explicit values of $n$ in the boundedness of $n$-complements for non-exceptional pairs of the same dimension.

\smallskip

\noindent\textbf{Shokurov's local index theorem}. An immediate application of Theorem \ref{thm: lc complement foliated surface local} is the local index theorem for foliated surfaces:
\begin{thm}\label{thm: surface index theorem local}
Let $a$ be a non-negative rational number and $\Ii\subset [0,1]\cap\mathbb Q$ a DCC set. Then there exists a positive integer $I$ depending only on $a$ and $\Ii$ satisfying the following. 

Assume that $(X\ni x,\Ff,B)$ is a foliated surface germ, such that $B\in\Ii$ and $\mld(X\ni x,\Ff,B)=a$. Then $I(K_{\Ff}+B)$ is Cartier near $x$. 
\end{thm}
Theorem \ref{thm: surface index theorem local} can be interpreted as a result in the direction of solving Shokurov's local index conjecture for foliated surfaces, which was posed in \cite[Conjecture 6.3]{CH21}. This conjecture is an important open problem in the study of foliations and has attracted significant attention in recent years. We provide a formal statement of the conjecture and its background in Conjecture \ref{conj: local index conjecture foliation} below.

\smallskip

\noindent\textbf{Minimal log discrepancies}. In this paper, we provide a characterization of the set of minimal log discrepancies of foliated surface singularities.
\begin{thm}\label{thm: set of foliated surface mlds}
Let $\Ii\subset [0,1]$ be a set. Then
\begin{align*}
    &\{\mld(X\ni x,\Ff,B)\mid \dim X=2,\rk\Ff=1,B\in\Ii\}\\
    =&\left\{0,\frac{1-\sum c_i\gamma_i}{n}\mid n\in\mathbb N^+,c_i\in\mathbb N,\gamma_i\in\Ii\right\}\cap [0,1].
\end{align*}
\end{thm}
Theorem \ref{thm: set of foliated surface mlds} implies the the following two results in \cite{Che23}.
\begin{cor}[$=${\cite[Remark 4.9]{Che23}}]\label{cor: set of foliated surface mlds no boundary}
$$\{\mld(X\ni x,\Ff)\mid \dim X=2,\rk\Ff=1\}=\left\{0,\frac{1}{n}\mid n\in\mathbb N^+\right\}.$$
\end{cor}
\begin{cor}[$=${\cite[Theorem 0.2]{Che23}}]\label{cor: acc mld foliation surface reprove}
Let $\Ii\subset [0,1]$ be a DCC set. Then
$$\{\mld(X\ni x,\Ff,B)\mid \dim X=2,\rk\Ff=1,B\in\Ii\}$$
satisfies the ACC.
\end{cor}

In addition to considering the possible values of mlds, it is also natural to investigate the structure of divisors that compute the mlds. We have the following result:

\begin{thm}\label{thm: uniform bound of mlds}
Let $\Ii\subset[0,1]$ a DCC set. Then there exists a positive real number $l$ depending only on $\Ii$ satisfying the following. 
    
Assume that $(X\ni x,\Ff,B)$ is an lc foliated surface germ such that $B\in\Ii$. Then there exists a prime divisor $E$ over $X\ni x$, such that $a(E,\Ff,B)=\mld(X\ni x,\Ff,B)$ and $a(E,\Ff,0)\leq l$.
\end{thm}
We remark that the condition ``$K_{\Ff}$ is $\Qq$-Cartier" is not necessary in Theorem \ref{thm: uniform bound of mlds}, as $\mld$ is well-defined for numerically lc foliations, as defined in Definition \ref{def: mld}.

Theorem \ref{thm: uniform bound of mlds} is a foliated surface case of the uniform boundedness conjecture for mlds, which can be found in \cite[Conjecture 8.2]{HLL22} (with an earlier form presented in \cite[Conjecture 1.1]{MN18}). More details on the conjecture and its background are provided in Conjecture \ref{conj: bdd mld computing divisor foliation} below.

\smallskip

\noindent\textbf{Uniform rational lc polytopes}. The last result in our paper is the existence of a uniform rational lc polytope for foliated surfaces. The theorem statement is as follows:

\begin{thm}\label{thm: uniform rational polytope foliation intro}
Let $v_1^0,\dots,v_m^0$ be positive integers and $\mathbf{v}_0:=(v_1^0,\dots,v_m^0)$. Then there exists an open set $U\ni\mathbf{v}_0$ of the rational envelope of $\mathbf{v}_0$ satisfying the following. 

Let $(X,\Ff,B=\sum_{i=1}^mv_i^0B_i)$ be any lc foliated triple of dimension $\leq 2$, where $B_i\geq 0$ are distinct Weil divisors. Then $(X,\Ff,B=\sum_{i=1}^mv_iB_i)$ is lc for any $(v_1,\dots,v_m)\in U$.
\end{thm}
Theorem \ref{thm: uniform rational polytope foliation intro} provides a positive answer to \cite[Conjecture 1.6]{LLM23} in dimension $2$, and it will be a key ingredient in our future work on the complete version (the real coefficients case) of the global ACC for foliated threefolds. Additionally, Theorem \ref{thm: uniform rational polytope foliation intro} can be viewed as the foliated surface case of the existence of uniform rational lc polytopes for usual pairs \cite[Theorem 5.6]{HLS19}.

\smallskip

\noindent\textit{Sketch of the paper}. In Section \ref{sec: Preliminaries} we introduce some preliminaries for foliations and also define complements for foliations. In Section \ref{sec: Preliminaries of foliations on surfaces} we recall the knowledge of foliations on surfaces, introduce and classify numerically lc foliated surface singularities. In Section \ref{sec: proof of main theorems} we prove all the other main theorems. In Section \ref{sec: foliation conjectures}, we formally state the foliated version of some standard conjectures in the minimal model program and discuss their background.

\section{Preliminaries}\label{sec: Preliminaries}

We work over the field of complex numbers $\mathbb C$. Our notation and definitions for algebraic geometry follow the standard references \cite{KM98, BCHM10}. For foliations, we adopt the notation and definitions introduced in \cite{LLM23}, which are based on those in \cite{CS20,ACSS21,CS21}.

\subsection{Foliations}

\begin{defn}[Foliations, {cf. \cite[Section 2.1]{CS21}}]\label{defn: foliation}
Let $X$ be a normal variety. A \emph{foliation} on $X$ is a coherent sheaf $\Ff\subset T_X$ such that
\begin{enumerate}
    \item $\Ff$ is saturated in $T_X$, i.e. $T_X/\Ff$ is torsion free, and
    \item $\Ff$ is closed under the Lie bracket.
\end{enumerate}
The \emph{rank} of the foliation $\Ff$ is the rank of $\Ff$ as a sheaf and is denoted by $\rk\Ff$. The \emph{co-rank} of $\Ff$ is $\dim X-\rk\Ff$. The \emph{canonical divisor} of $\Ff$ is a divisor $K_\Ff$ such that $\mathcal{O}_X(-K_{\mathcal{F}})\cong\mathrm{det}(\Ff)$. 
\end{defn}

\begin{defn}[Singular locus]
     Let $X$ be a normal variety and let $\Ff$ be a rank $r$ foliation on $X$. We can associate to $\Ff$ a morphism $$\phi: \Omega_X^{[r]}\to \mathcal{O}_X(K_{\Ff})$$ defined by taking the double dual of the $r$-wedge product of the map $\Omega^1_X\to \Ff^*$, induced by the inclusion $\Ff\to T_X$. This yields a map $$\phi': (\Omega_X^{[r]}\otimes\mathcal{O_X}(-K_{\Ff}))^{\vee\vee}\to \mathcal{O}_X$$ and we define the singular locus, denoted as $\Sing~\Ff$, to be the co-support of the image of $\phi'$.
\end{defn}

\begin{defn}[Pullback and pushforward, {cf. \cite[3.1]{ACSS21}}]\label{defn: pullback}
Let $X$ be a normal variety, $\Ff$ a foliation on $X$, $f: Y\dashrightarrow X$ a dominant map, and $g: X\dashrightarrow X'$ a birational map. We denote $f^{-1}\Ff$ the \emph{pullback} of $\Ff$ on $Y$ as constructed in \cite[3.2]{Dru21}. We also say that $f^{-1}\Ff$ is the \emph{induced foliation} of $\Ff$ on $Y$. We define $g_*\Ff:=(g^{-1})^{-1}\Ff$ and denote it by $g_*\Ff$.
\end{defn}

\begin{defn}[Invariant subvarieties, {cf. \cite[3.1]{ACSS21}}]\label{defn: f-invariant}
Let $X$ be a normal variety, $\Ff$ a foliation on $X$, and $S\subset X$ a subvariety. We say that $S$ is \emph{$\Ff$-invariant} if and only if for any open subset $U\subset X$ and any section $\partial\in H^0(U,\Ff)$, we have $$\partial(\mathcal{I}_{S\cap U})\subset \mathcal{I}_{S\cap U}$$ 
where $\mathcal{I}_{S\cap U}$ is the ideal sheaf of $S\cap U$. 
\end{defn}

\begin{defn}[Non-dicritical singularities, {cf. \cite[Definition 2.10]{CS21}}]\label{defn: non-dicritical}
Let $X$ be a normal variety and $\Ff$ a foliation of co-rank $1$ on $X$. We say that $\Ff$ has \emph{non-dicritical} singularities if for any closed point $x\in X$ and any birational morphism $f: X'\rightarrow X$ such that $f^{-1}(\overline{\{x\}})$ is a divisor, each component of $f^{-1}(\overline{\{x\}})$ is $f^{-1}\Ff$-invariant.
\end{defn}

\begin{defn}[Special divisors on foliations, cf. {\cite[Definition 2.2]{CS21}}]\label{defn: special divisors on foliations}
Let $X$ be a normal variety and $\Ff$ a foliation on $X$. For any prime divisor $C$ on $X$, we define $\epsilon_{\Ff}(C):=1$ if $C$ is not $\Ff$-invariant, and  $\epsilon_{\Ff}(C):=0$ if $C$ is $\Ff$-invariant. If $\Ff$ is clear from the context, then we may use $\epsilon(C)$ instead of $\epsilon_{\Ff}(C)$. For any $\Rr$-divisor $D$ on $X$, we define $$D^{\Ff}:=\sum_{C\text{ is a component of }D}\epsilon_{\Ff}(C)\cdot C.$$
Let $E$ be a prime divisor over $X$ and $f: Y\rightarrow X$ a projective birational morphism such that $E$ on $Y$. We define $\epsilon_{\Ff}(E):=\epsilon_{f^{-1}\Ff}(E)$. It is clear that $\epsilon_{\Ff}(E)$ is independent of the choice of $f$.
\end{defn}

\begin{defn}
A \emph{foliated sub-triple} $(X/Z\ni z,\Ff,B)$ consists of a normal quasi-projective variety $X$, a foliation $\Ff$ on $X$, an $\Rr$-divisor $B$ on $X$, and a projective morphism $X\rightarrow Z$, and a (not necessarily closed) point $z\in Z$, such that $K_{\Ff}+B$ is $\Rr$-Cartier over a neighborhood of $z$. 

If $\Ff=T_X$, then we may drop $\Ff$ and denote $(X/Z\ni z,\Ff,B)$ by $(X/Z\ni z,B)$, and say that $(X/Z\ni z,B)$ is a \emph{sub-pair}. If $B\geq 0$ over a neighborhood of $z$, then we say that $(X/Z\ni z,\Ff,B)$ is a \emph{foliated triple}. If $\Ff=T_X$ and  $B\geq 0$ over a neighborhood of $z$, then we say that $(X/Z\ni z,B)$ is a \emph{pair}.

Let $(X/Z\ni z,\Ff,B)$ be a foliated (sub-)triple. If $X\rightarrow Z$ is the identity morphism, then we may drop $Z$ and denote $(X/Z\ni z,\Ff,B)$ by $(X\ni z,\Ff,B)$, and say that $(X\ni z,B)$ is a \emph{foliated (sub-)germ}. If $\Ff=T_X$ and $X\rightarrow Z$ is the identity morphism,  we may drop $Z$ and say that $(X\ni z,B)$ is a \emph{(sub-)germ}. 

If $(X/Z\ni z,\Ff,B)$ (resp. $(X\ni z,\Ff,B)$, $(X/Z\ni z,B)$, $(X\ni z,B)$) is an foliated (sub-)triple (resp. foliated (sub-)germ, (sub-)pair, (sub-)germ) for any $z\in Z$, then we say that  $(X/Z,\Ff,B)$ (resp. $(X,\Ff,B)$, $(X/Z,B)$, $(X,B)$) is a foliated (sub-)triple (resp. foliated (sub-)triple, (sub-)pair, (sub-)pair).
\end{defn}

\begin{defn}\label{defn: sing of f triples}
Let $(X/Z\ni z,\Ff,B)$ be a foliated (sub-)triple.  For any prime divisor $E$ over $X$, let $f: Y\rightarrow X$ be a birational morphism such that $E$ is on $Y$, and suppose that
$$K_{\Ff_Y}+B_Y=f^*(K_\Ff+B)$$
over a neighborhood of $z$, where $\Ff_Y:=f^{-1}\Ff$. We define $a(E,\Ff,B):=-\mult_EB_Y$ to be the \emph{discrepancy} of $E$ with respect to $(X,\Ff,B)$. It is clear that $a(E,\Ff,B)$ is independent of the choice of $Y$. If $\Ff=T_X$, then we let $a(E,X,B):=a(E,\Ff,B)$.

Let $\delta$ be a non-negative real number and $(X/Z\ni z,\Ff,B)$ a foliated (sub-)triple, We say that $(X/Z\ni z,\Ff,B)$ is \emph{(sub-)lc} (resp. \emph{(sub-)klt}, \emph{(sub-)$\delta$-lc}, \emph{(sub-)$\delta$-klt}, \emph{(sub-)canonical}, \emph{(sub-)terminal}) if $a(E,\Ff,B)\geq -\epsilon_{\Ff}(E)$ (resp. $>-\epsilon_{\Ff}(E)$, $\geq-\epsilon_{\Ff}(E)+\delta$, $>-\epsilon_{\Ff}(E)+\delta$, $\geq 0$, $>0$) for any prime divisor $E$ over $z$, i.e. the closure of the image of $E$ on $Z$ is $\bar z$. We define
$$\mld(X/Z\ni z,\Ff,B):=\inf\{a(E,\Ff,B)+\epsilon_{\Ff}(E)\mid E\text{ is over }z\}$$
to be the \emph{minimal log discrepancy} (\emph{mld} for short) of $\mld(X/Z\ni z,\Ff,B)$.

Let $(X,\Ff,B)$ be a foliated (sub-)triple. We say that $(X,\Ff,B)$ is \emph{(sub-)lc} (resp. \emph{(sub-)klt}, \emph{(sub-)$\delta$-lc}, \emph{(sub-)$\delta$-klt}) if $(X\ni x,\Ff,B)$ is \emph{(sub-)lc} (resp. \emph{(sub-)klt}, \emph{(sub-)$\delta$-lc}, \emph{(sub-)$\delta$-klt}) for any point $x\in X$. We say that  $(X,\Ff,B)$ is \emph{(sub-)canonical} (resp. \emph{(sub-)terminal}) if $(X\ni x,\Ff,B)$ is \emph{(sub-)canonical} (resp. \emph{(sub-)terminal}) for any codimension $\geq 2$ point $x\in X$.
\end{defn}

\subsection{Complements}\label{section: foliation complements}
\begin{defn}\label{defn: complement}
Let $n$ be a positive integer, $\epsilon$ a non-negative real number, $\Ii_0\subset (0,1]$ a finite set, and $(X/Z\ni z,\Ff,B)$ and $(X/Z\ni z,\Ff,B^+)$ two foliated triples. We say that $(X/Z\ni z,\Ff,B^+)$ is an \emph{$(\epsilon,\Rr)$-complement} of $(X/Z\ni z,\Ff,B)$ if 
\begin{itemize}
    \item $(X/Z\ni z,\Ff,B^+)$ is $\epsilon$-lc,
    \item $B^+\geq B$, and
    \item $K_{\Ff}+B^+\sim_{\Rr}0$ over a neighborhood of $z$.
\end{itemize}
We say that $(X/Z\ni z,\Ff,B^+)$ is an \emph{$(\epsilon,n)$-complement} of $(X/Z\ni z,\Ff,B)$ if
\begin{itemize}
\item $(X/Z\ni z,\Ff,B^+)$ is $\epsilon$-lc,
\item $nB^+\geq \lfloor (n+1)\{B\}\rfloor+n\lfloor B\rfloor$, and
\item $n(K_{\Ff}+B^+)\sim 0$ over a neighborhood of $z$.
\end{itemize}
 We say that $(X/Z\ni z,\Ff,B)$ is $(\epsilon,\Rr)$-complementary if $(X/Z\ni z,\Ff,B)$ has an $(\epsilon,\Rr)$-complement. We say that $(X/Z\ni z,\Ff,B^+)$ is a \emph{monotonic $(\epsilon,n)$-complement} of $(X/Z\ni z,\Ff,B)$ if $(X/Z\ni z,\Ff,B^+)$ is an $(\epsilon,n)$-complement of $(X/Z\ni z,\Ff,B)$ and $B^+\geq B$.
 
$(0,\Rr)$-complement (resp. $(0,n)$-complement, $(0,\Rr)$-complementary, $(0,n)$-complementary) is also called $\Rr$-complement (resp. $n$-complement, $\Rr$-complementary, $n$-complementary). 
\end{defn}

\section{Foliations on surfaces}\label{sec: Preliminaries of foliations on surfaces}

\subsection{Resolution of foliated surfaces}

\begin{defn}
Let $X$ be a normal surface, $\Ff$ a foliation on $X$, and $x\in X$ a closed point such that $x\not\in\Sing(X)$ and $x\in\Sing(\Ff)$. Let $v$ be a vector field generating $\Ff$ near $x$. By \cite[Page 2, Line 17-18]{Bru15}, $v(x)=0$ and $(Dv)|_x$ has exactly two eigenvalues $\lambda_1$ and $\lambda_2$. 

We say that $x$ is a \emph{reduced singularity} of $\Ff$ if at least one of $\lambda_1$ and $\lambda_2$ is not $0$ (say, $\lambda_2$) and $\frac{\lambda_1}{\lambda_2}\not\in\mathbb Q^+$. We say that $x$ is a \emph{non-degenerate reduced singularity} of $\Ff$ if $x$ is a reduced singularity of $\Ff$ and $\frac{\lambda_1}{\lambda_2}\not\in\{0,\infty\}$, i.e. $\lambda_1$ and $\lambda_2$ are both not equal to $0$.

We say that $\Ff$ has \emph{at most reduced singularities} if for any closed point $p\in X$, $\Ff$ is either non-singular at $p$ or $p$ is a reduced singularity of $\Ff$. 

An \emph{$\Ff$-exceptional curve} is a non-singular rational curve $E$ on $X$ such that 
\begin{enumerate}
    \item $X$ is smooth near $E$ and $E^2=-1$, 
    \item there exists a divisorial contraction $f: X\rightarrow Y$ of $E$, and
    \item  $f(E)$ is a reduced singularity of $f_*\Ff$.
\end{enumerate}
\end{defn}

\begin{defn}[Minimal resolution]\label{defn: minimal resolution foliation}
Let $X$ be a normal surface, $\Ff$ a foliation on $X$, $f: Y\rightarrow X$ a projective birational morphism, and $\Ff_Y:=f^{-1}\Ff$.

We say that $f$ is a \emph{resolution} of $\Ff$ if $Y$ is smooth and $\Ff_Y$ has at most reduced singularities. We say that  $f$ is the \emph{minimal resolution} of $\Ff$ if for any resolution $g: W\rightarrow X$ of $\Ff$, $g$ factors through $f$, i.e. there exists a projective birational morphism $h: W\rightarrow Y$ such that $g=f\circ h$. By definition, the minimal resolution of $\Ff$ is unique if it exists.

%We say that $f$ is a \emph{proper resolution} of $\Ff$ if $Y$ is smooth and $\Ff_Y$ has at most reduced singularities. We say that  $f$ is a \emph{minimal resolution} of $\Ff$ if for any proper resolution $g: W\rightarrow X$ of $\Ff$, $g$ factors through $f$, i.e. there exists a projective birational morphism $h: W\rightarrow Y$ such that $g=f\circ h$. By definition, if $\Ff$ has a minimal resolution, then it is unique, so we may call a minimal resolution of $\Ff$ as \emph{the} minimal resolution of $\Ff$.

For any closed point $x\in X$, the \emph{minimal resolution} of $\Ff\ni x$ is the minimal resolution of $\Ff$ for any sufficiently small neighborhood of $x$.
\end{defn}

\begin{prop}[{\cite[Proposition 1.17]{Che23}}]\label{prop: existence minimal resolution}
Let $X$ be a normal surface and $\Ff$ a foliation on $X$. Then the minimal resolution of $\Ff$ exists.

%Then $\Ff$ has a unique minimal resolution.
\end{prop}

\subsection{Invariants of curves on foliated surfaces}

\begin{defn}
Let $X$ be a normal surface with at most cyclic quotient singularities, $\Ff$ a foliation on $X$, and $C$ a reduced curve on $X$ such that no component of $C$ is $\Ff$-invariant. For any closed point $x\in X$, we define $\tang(\Ff,C,x)$ in the following way. 
\begin{itemize}
\item If $x\notin \Sing(X)$, then we let $v$ be a vector field generating $\Ff$ around $x$, and $f$ a holomorphic function defining $C$ around $x$. We define 
$$\tang(\Ff,C,x):=\dim_{\mathbb{C}}\frac{\mathcal{O}_{X,x}}{\langle f, v(f)\rangle}.$$
\item If $x\in\Sing(X)$, then $x$ is a cyclic quotient singularity of index $r$ for some integer $r\geq 2$. Let $\rho:\tilde X\rightarrow X$ be an index $1$ cover of $X\ni x$, $\tilde x:=\rho^{-1}(x)$, $\widetilde C:=\rho^*C$, and $\tilde\Ff$ the foliation induced by the sheaf $\rho^*\Ff$ near $\tilde x$. Then $\tilde x$ is a smooth point of $\tilde X$, and we define
$$\tang(\Ff,C,x):=\frac{1}{r}\tang(\tilde\Ff,\tilde C,\tilde x).$$
\end{itemize}
We define
$$\tang(\Ff,C):=\sum_{x\in X}\tang(\Ff,C,x).$$
By \cite[Section 2]{Bru02}, $\tang(\Ff,C)$ is well-defined.
\end{defn}

\begin{defn}
%Let $X$ be a normal surface with at most cyclic quotient singularities, 
Let $X$ be a smooth surface,
$\Ff$ a foliation on $X$, and $C$ a reduced curve on $X$ such that all components of $C$ are $\Ff$-invariant. For any closed point $x\in X$, we define $Z(\Ff,C,x)$ in the following way.

Let $\omega$ be a $1$-form generating $\Ff$ around $x$, and $f$ a holomorphic function generating $C$ around $x$. Then there are uniquely determined holomorphic functions $g,h$ and a holomorphic $1$-form $\eta$ on $X$ near $x$, such that $g\omega=hdf+f\eta$ and $f,h$ are coprime.
We define 
$$Z(\Ff,C,x):=\text{ the vanishing order of } \frac{h}{g}|_C \text{ at } x.$$
By \cite[Chapter 2, Page 15]{Bru15}, $Z(\Ff,C,x)$ is independent of the choice of $\omega$.

%\begin{itemize}
%\item If $x\notin \Sing(X)$, then we let $\omega$ be a $1$-form generating $\Ff$ around $x$, and $f$ a holomorphic function generating $C$ around $x$. Then there are uniquely determined holomorphic functions $g,h$ and a holomorphic $1$-form $\eta$ on $X$ near $x$, such that $g\omega=hdf+f\eta$ and $f,h$ are coprime.
%We define 
%$$Z(\Ff,C,x):=\text{ the vanishing order of } \frac{h}{g}|_C \text{ at } x.$$
%By \cite[Chapter 2, Page 15]{Bru15}, $Z(\Ff,C,x)$ is independent of the choice of $\omega$.
%\item If $x\in C\cap \Sing(X)$, we define 
%$Z(\Ff,C,x):=0.$
%\end{itemize}
We define $$Z(\Ff,C):=\sum_{x\in C}Z(\Ff,C,x).$$
%By \cite[Section 2]{Bru02}, $Z(\Ff,C)$ is well-defined.
\end{defn}

\begin{thm}[{cf. \cite[Chapter 2]{Bru15}}]\label{thm: tang and Z formula}
Let $X$ be a smooth quasi-projective surface, %with at most cyclic quotient singularities,
$\Ff$ a foliation on $X$, and $C$ a compact reduced curve on $X$. %Suppose that $\Sing(\Ff)\cap\Sing(X)\cap C=\emptyset$.
\begin{enumerate}
    \item If no components of $C$ is $\Ff$-invariant, then 
    $$K_{\Ff}\cdot C+C^2=\tang(\Ff,C).$$
    \item If all components of $C$ are $\Ff$-invariant, then
    $$K_{\Ff}\cdot C=Z(\Ff,C)-\chi(C)$$
    where $\chi(C):=-K_X\cdot C-C^2$.
\end{enumerate}
\end{thm}

The following lemma is a variation of Theorem \ref{thm: tang and Z formula}(1).

\begin{lem}[{cf. \cite[Proposition 2.2]{Bru15}}]\label{lem: restrction tang}
Let $X$ be a smooth surface, $\Ff$ a foliation on $X$, and $C$ a compact reduced curve on $X$ such that no component of $C$ is $\Ff$-invariant. Then there exists a Weil divisor $D\geq 0$ on $C$ such that $(K_{\Ff}+C)|_C\sim D$ and $\deg D=\tang(\Ff,C)$.
\end{lem}
\begin{proof}
We choose an open covering $\{U_j\}$ of $X$ , holomorphic vector field $v_j$ on $U_j$ generating $\Ff$, and holomorphic functions $f_j$ on $U_j$ defining $C$. On the intersections $U_i\cap U_j$ we have $v_i=g_{i,j}v_j$ and $f_i=f_{i,j}f_j$, where $g_{i,j}$ are cocycles representing $T_{\Ff}^\vee\cong\mathcal{O}_X(K_{\Ff})$ and $f_{i,j}$ are cocycles representing $\mathcal{O}_X(C)$. Hence the functions $\{v_j(f_j)\}$ restricted to $C$ give a section of $(T_{\Ff}^\vee\otimes\mathcal{O}_X(C))|_C$, because by Leibniz’s rule,
$$v_i(f_i)=g_{i,j}v_j(f_{i,j}f_j)=g_{i,j}f_{i,j}v_j(f_j)+g_{i,j}f_jv_j(f_{i,j})$$
and $g_{i,j}f_jv_j(f_{i,j})=0$ on $C$. We let $D$ be a section of $(T_{\Ff}^\vee\otimes\mathcal{O}_X(C))|_C$, then $(K_{\Ff}+C)|_C\sim D\geq 0$. Moreover, $D$ vanishes at the points of $C$ where $\Ff$ is not transverse to $C$, and the vanishing order is nothing but that $\tang(\Ff,C)$. In other words, $\deg D=\tang(\Ff,C)$.
\end{proof}

\subsection{Dual graphs}

\begin{defn}[Dual graph]\label{defn: dual graph}
Let $n$ be a non-negative integer, and $C=\cup_{i=1}^nC_i$ a collection of irreducible curves contained in the non-singular locus of a normal surface $X$. We define the \emph{dual graph} $\mathcal{D}(C)$ of $C$ as follows.
\begin{enumerate}
    \item The vertices $v_i=v_i(C_i)$ of $\mathcal{D}(C)$ correspond to the curves $C_i$.
    \item For $i\neq j$, the vertices $v_i$ and $v_j$ are connected by $C_i\cdot C_j$ edges.
    \item Each vertex $v_i$ is labeled by $w(C_i):=-C_i^2$. The integer $w(C_i)$ is called the \emph{weight} of $C_i$.
\end{enumerate}
We sometimes write the name of the curve $C_i$ near the vertex $v_i$. We say that $\mathcal{D}(C)$ is a \emph{cycle} if
\begin{itemize}
    \item either $n=2$ and $C_1\cdot C_2=2$, or
    \item $n\geq 3$, and possibly reordering indices, we have 
    \begin{itemize}
        \item $C_i\cdot C_j=1$ if  $|i-j|=1$ or $\{i,j\}=\{1,n\}$, and
        \item $C_i\cdot C_j=0$ if $|i-j|\geq 2$ and $\{i,j\}\not=\{1,n\}$.
    \end{itemize}
\end{itemize}
We say that $\mathcal{D}(C)$ contains a \emph{cycle} if there exists a sub-dual graph of $\mathcal{D}(C)$ that is a cycle. We say that $\mathcal{D}(C)$ is a \emph{tree} if
\begin{itemize}
    \item $\mathcal{D}(C)$ does not contain a cycle, and
    \item $C_i\cdot C_j\leq 1$ for any $i\not=j$.
\end{itemize}
The \emph{intersection matrix} of $\mathcal{D}(C)$ is defined as the matrix $(C_i\cdot C_j)_{1\leq i,j\leq n}$ if $C\not=\emptyset$. The \emph{determinant} of $\mathcal{D}(C)$ is defined as
$$\det(\mathcal{D}(C)):=\det(-(C_i\cdot C_j)_{1\leq i,j\leq n})$$
if $C\not=\emptyset$, and $\det(\mathcal{D}(C)):=1$ if $C=\emptyset$.

A \emph{fork} of $\mathcal{D}(C)$ is a curve $C_i$ such that $C_i\cdot C_j\geq 1$ for at least three different $j\not=i$, and we also say that $v_i$ is a \emph{fork}. A \emph{tail} of $\mathcal{D}(C)$ is a curve $C_i$ such that $C_i\cdot C_j\geq 1$ for at most one $j\not=i$, and we also say that $v_i$ is a \emph{tail}. A \emph{chain} is a dual graph that is a tree which does not contain a fork.

For any $i,j$, we say that $C_i$ and $C_j$ are \emph{adjacent} if $i\not=j$ and $C_i\cdot C_j\geq 1$.

For any projective birational morphism $f: Y\rightarrow X$ between surfaces, let $E=\cup_{i=1}^nE_i$ be the reduced exceptional divisor for some non-negative integer $n$. Suppose that $E$ is contained in the non-singular locus of $Y$. Then we define $\mathcal{D}(f):=\mathcal{D}(E)$.
\end{defn}

\begin{defn}[Dual graph on a foliated surface]
Let $n$ be a positive integer, $X$ a normal surface, $\Ff$ a foliation on $X$, and $C=\cup_{i=1}^nC_i$ a collection of irreducible curves contained in the non-singular locus of $X$. 
\begin{enumerate}
    \item We say that $C=\cup_{i=1}^nC_i$ is a \emph{string} if
    \begin{enumerate}
        \item for any $i$, $C_i$ is a smooth rational curve, and
        \item for any $i,j$, 
        \begin{enumerate}
            \item $C_i\cdot C_j=1$ if $|i-j|=1$, and
            \item $C_i\cdot C_j=0$ if $|i-j|>1$.
        \end{enumerate}
    \end{enumerate}
    \item We say that $C=\cup_{i=1}^nC_i$ is a \emph{Hirzebruch-Jung string} if $C=\cup_{i=1}^nC_i$ is a string and $C_i^2\leq -2$ for any $i$.
    \item We say that $C=\cup_{i=1}^nC_i$ is an $\Ff$-chain if
    \begin{enumerate}
        \item $C=\cup_{i=1}^nC_i$ is a Hirzebruch-Jung string,
        \item $C_i$ is $\Ff$-invariant for any $i$,
        \item for any closed point $x\in C$, either $x\not\in\Sing(\Ff)$, or $x$ is a non-degenerate reduced singularity of $\Ff$, and
        \item $Z(\Ff,C_1)=1$, and $Z(\Ff,C_i)=2$ for any $i\geq 2$.
    \end{enumerate}
\end{enumerate}
\end{defn}

\subsection{Surface foliated numerical triples}

\cite[Fact I.2.4]{McQ08} have classified all foliated surface singularities $(X\ni x,\Ff)$ such that $\Ff$ is canonical at $x$, while \cite[Theorem 0.1]{Che23} has classified all foliated surface singularities $(X\ni x,\Ff)$ such that $\Ff$ is lc at $x$. However, in practice, we may usually come up with the structure of lc foliated germs $(X\ni x,\Ff,B)$. If $B\not=0$, then we don't know whether $K_{\Ff}$ is $\Qq$-Cartier near $x$ or not, so the forementioned classification results cannot be directly applied. To resolve this issue, we need to introduce the concept of numerical surface singularities of foliations and provide the classification of numerically lc surface singularities of foliations.

From now until the end of this section, we present a detailed characterization of lc foliated surface singularities. We define numerically lc (num-lc for short) foliated surface singularities similar to the definition of num-lc singularities for usual surface singularities. We then classify all num-lc foliated surface singularities. Although the result is very similar to \cite[Theorem 0.1]{Che23}, for the reader's convenience, we provide a complete and detailed proof in this section.

\begin{defn}\label{defn: num foliated sing}
A \emph{surface foliated numerical sub-triple} (\emph{surface foliated num-sub-triple} for short) $(X,\Ff,B)$ consists of a normal surface $X$, a rank $1$ foliation $\Ff$ on $X$, and an $\Rr$-divisor $B$ on $X$. We say that $(X,\Ff,B)$ is a \emph{surface foliated numerical triple} (\emph{surface foliated num-triple} for short) if  $(X,\Ff,B)$ is a surface foliated num-sub-triple and $B\geq 0$. A \emph{surface foliated numerical germ} (\emph{surface foliated num-germ} for short) $(X\ni x,\Ff,B)$ consists of a surface foliated num-triple $(X,\Ff,B)$ and a closed point $x\in X$.

Let $(X,\Ff,B)$ be a surface foliated num-sub-triple. Let $f: Y\rightarrow X$ of $X$ be a resolution of $X$ with prime $f$-exceptional divisors $E_1,\dots,E_n$ for some non-negative integer $n$ such that $\Center_YE$ is a divisor. Since  $\{(E_i\cdot E_j)\}_{n\times n}$ is negative definite, the equation 
$$\begin{pmatrix}
  (E_1\cdot E_1) &\cdots & (E_1\cdot E_n) \\
  \vdots&\ddots & \vdots         \\
  (E_n\cdot E_1) &\cdots & (E_n\cdot E_n) \\
\end{pmatrix} 
\begin{pmatrix}
  a_1 \\
  \vdots \\
  a_n \\
\end{pmatrix} 
 =
 \begin{pmatrix}
-(K_{\Ff_Y}+B_Y)\cdot E_1\\
\vdots\\
-(K_{\Ff_Y}+B_Y)\cdot E_n\\
\end{pmatrix}$$
has a unique solution $(a_1,\dots,a_n)$, where $\Ff_Y:=f^{-1}\Ff$ and $B_Y:=f^{-1}_*B$. We define $$a_{\num,f}(E,\Ff,B):=-\mult_E\left(B_Y+\sum_{i=1}^na_iE_i\right).$$
\end{defn}

As in the pair case \cite[Chapter 4]{KM98}, it is easy to see that $a_{\num,f}(E,\Ff,B)=a(E,\Ff,B)$ when $K_\Ff+B$ is $\Rr$-Cartier and $a_{\num,f}(E,\Ff,B)$ does not depend on the resolution $f$. This enable us to define $a(E,\Ff,B)$ when $K_\Ff+B$ is not necessarily $\Rr$-Cartier.

\vspace{.5em}

\begin{defn}
Let $(X,\Ff,B)$ be a surface foliated num-sub-triple. We define $a(E,\Ff,B):=a_{\num,f}(E,\Ff,B)$ for an arbitrary resolution $f: Y\rightarrow X$ of $X$ such that $E$ is a divisor on $Y$. %Lemmas \ref{lem: anum same as a} and \ref{lem: anum not depend on resolution} guarantee that there is no abuse of notations.

Let $(X\ni x,\Ff,B)$ be a surface foliated num-germ. We say that $(X\ni x,\Ff,B)$ is \emph{num-lc} (resp. \emph{num-klt, num-canonical, num-terminal}) if $a(E,\Ff,B)\geq-\epsilon(E)$ (resp. $>-\epsilon(E),\geq 0,>0$) for any prime divisor $E$ over $X\ni x$. If $K_\Ff+B$ is $\Rr$-Cartier at $x$, then num-lc (resp. num-klt, num-canonical, num-terminal) is equivalent to lc (resp. klt, canonical, terminal).
\end{defn}

\begin{comment}
\begin{lem}\label{lem: num lc and rcartier imply lc}
Let $(X\ni x,\Ff,B)$ be an foliated germ such that $\dim X=2$, $\rk\Ff=1$, and $x$ is a closed point. If $(X\ni x,\Ff,B)$ is \emph{num-lc} (resp. \emph{num-klt, num-canonical, num-terminal}), then $(X\ni x,\Ff,B)$ is \emph{lc} (resp. \emph{klt, canonical, terminal}).
\end{lem}
\begin{proof}
Notice that $K_\Ff+B$ is $\Rr$-Cartier at $x$ by definition. The lemma immediately follows from Lemma \ref{lem: anum same as a}.
\end{proof}
\end{comment}

\begin{defn}\label{def: mld}
Let $(X\ni x,\Ff,B)$ be a surface foliated num-germ. The \emph{minimal log discrepancy} (\emph{mld} for short) of $(X\ni x,\Ff,B)$ is defined as $$\mld(X\ni x,\Ff,B):=\inf\{a(E,\Ff,B)+\epsilon_{\Ff}(E)\mid E \text{ is a prime divisor over }X\ni x\}.$$
%By Lemma \ref{lem: anum not depend on resolution}, 
We define $\mld(X\ni x,\Ff):=\mld(X\ni x,\Ff,0)$. Notice that this definition coincides with the mld defined in Definition \ref{defn: sing of f triples} for foliated germs. 
\end{defn}

\begin{lem}
Let $(X\ni x,\Ff,B)$ be a surface foliated num-germ. Then either $\mld(X\ni x,\Ff,B)=-\infty$, or 
$$\mld(X\ni x,\Ff,B)=\min\{a(E,\Ff,B)+\epsilon_{\Ff}(E)\mid E \text{ is a prime divisor over }X\ni x\}\geq 0.$$
\end{lem}
\begin{proof}
First suppose that $\mld(X\ni x,\Ff,B)<0$. Then there exists a resolution $f: Y\rightarrow X$ of $X\ni x$ with prime $f$-exceptional divisors $E_1,\dots,E_n$, and a prime divisor $E$ on $Y$, such that $\Center_XE=x$ and $\mult_EB_Y>\epsilon_{\Ff}(E)$, where $B_Y:=f^{-1}_*B-\sum_{i=1}^na(E_i,\Ff,B)\cdot E_i$. Then $(Y,\Ff_Y,B_Y)$ is not sub-lc near $E$ and by \cite[Proposition 3.4]{Che23} $\mld(Y,\Ff_Y,B_Y)=-\infty$, hence we also have $\mld(X\ni x,\Ff,B)=-\infty$.
%so for any positive integer $m$, there exists a prime divisor $F_m$ over $Y$ such that $a(F_m,\Ff_Y,B_Y)<-m$ and $\Center_{Y}F_m\in E$. Thus $\Center_XF_m=x$ for any $m$, so $\mld(X\ni x,\Ff,B)=-\infty$.

Now we suppose that  $\mld(X\ni x,\Ff,B)\geq 0$. Let $f: Y\rightarrow X$ of $X\ni x$ be a resolution of $X$ with prime $f$-exceptional divisors $E_1,\dots,E_n$, and let $B_Y:=f^{-1}_*B-\sum_{i=1}^na(E_i,\Ff,B)\cdot E_i$. Since $(Y,\Ff_Y,B_Y)$ is a foliated sub-triple over a neighborhood of $x$, then by \cite[Corollary 3.6]{Che23} we have 
\begin{align*}
  \mld(X\ni x,\Ff,B)&=\inf\{a(E,\Ff_Y,B_Y)+\epsilon_{\Ff}(E)\mid E \text{ is a prime divisor over }X\ni x\}\\
  &=\min\{a(E,\Ff_Y,B_Y)+\epsilon_{\Ff}(E)\mid E \text{ is a prime divisor over }X\ni x\}\\
  &=\min\{a(E,\Ff,B)+\epsilon_{\Ff}(E)\mid E \text{ is a prime divisor over }X\ni x\}.
\end{align*}
\end{proof}

\begin{lem}\label{lem: add component worse sing}
Let $(X\ni x,\Ff,B)$ be a surface num-lc foliated num-germ. Suppose that all components of $B$ pass through $x$. Then $\mld(X\ni x,\Ff,B)\leq\mld(X\ni x,\Ff,0)$, and $\mld(X\ni x,\Ff,B)<\mld(X\ni x,\Ff)$ if $B\not=0$. In particular, $(X\ni x,\Ff)$ is num-lc.
\end{lem}
\begin{proof} By Definition \ref{defn: num foliated sing} and \cite[Lemma 3.41]{KM98}, $a(E,\Ff,B)\le a(E,\Ff,0)$ for any prime divisor $E$ over $X\ni x$, and $a(E,\Ff,B)< a(E,\Ff,0)$ for any prime divisor $E$ over $X\ni x$ if $B\not=0$.
\end{proof}

\begin{defn}\label{def: pld}
Let $(X\ni x,\Ff,B)$ be a surface foliated num-germ and $f: Y\rightarrow X$ the minimal resolution of $\Ff\ni x$ such that $f$ is not the identity morphism. The \emph{partial log discrepancy} (\emph{pld} for short) of $(X\ni x,\Ff,B)$ is defined as
$$\pld(X\ni x,\Ff,B):=\min\{a(E,\Ff,B)+\epsilon_{\Ff}(E)\mid E \text{ is a }f\text{-exceptional prime divisor}\}.$$
We define $\pld(X\ni x,\Ff):=\pld(X\ni x,\Ff,0)$.
\end{defn}

Finally, we are ready to state and prove the main theorem of this section, which is a generalization of the classification theorems in \cite{McQ08} and \cite{Che23}, but with more details:

\begin{thm}\label{thm: classification of foliation intro}
Let $(X\ni x,\Ff,B)$ be a numerically lc surface foliated numerical germ such that all components of $B$ pass through $x$ and $\rk\Ff=1$. Let $f: Y\rightarrow X$ be the minimal resolution of $\Ff\ni x$ (cf. Definition \ref{defn: minimal resolution foliation}), $\mathcal{D}$ the dual graph of $f$, and $\Ff_Y:=f^{-1}\Ff$. Suppose that $f$ is not the identity morphism. Then one of the following cases holds.
\begin{enumerate}
    \item [(Case 1)] $\mathcal{D}=\cup_{i=1}^mE_i$ is an $\Ff_Y$-chain. Moreover, in this case,
    \begin{enumerate}
        \item $X\ni x$ is a cyclic quotient singularity and $\Ff$ is non-dicritical near $x$. In particular, $X\ni x$ is klt and $K_{\Ff}$ is $\Qq$-Cartier near $x$,
        \item for any $1\leq i\leq m$,
        $$a(E_i,\Ff)=\frac{\det(\mathcal{D}(\cup_{j=i+1}^mE_i))}{\det(\mathcal{D})},$$
    
        \item $\pld(X\ni x,\Ff)=a(E_m,\Ff)=\frac{1}{\det(\mathcal{D})}>0$, and $(X\ni x,\Ff)$ is terminal, and
        \item there is a unique $\Ff$-invariant curve $C$ passing through $x$ and $C$ is smooth at $x$.
        %\item $\Ff$ is non-singular near $x$.
    \end{enumerate}
    \item [(Case 2)] $\mathcal{D}=\cup_{i=1}^3E_i$ is a Hirzebruch-Jung string such that $Z(\Ff_Y,E_1)=Z(\Ff_Y,E_3)=1,Z(\Ff_Y,E_2)=3$, $E_1^2=E_3^2=-2$, and $E_2^2\leq -2$. Moreover, in this case,
    \begin{enumerate}
        \item $X\ni x$ is a cyclic quotient singularity and $\Ff$ is non-dicritical near $x$. In particular, $X\ni x$ is klt and $K_{\Ff}$ is $\Qq$-Cartier near $x$,
        \item 
        $a(E_1,\Ff)=a(E_3,\Ff)=\frac{1}{2}$ and $a(E_2,\Ff)=0$,
        \item $\pld(X\ni x,\Ff)=0$, $B=0$, and $(X\ni x,\Ff)$ is canonical but not terminal,and
        \item $2K_{\Ff}$ is Cartier near $x$.
    \end{enumerate}
    \item [(Case 3)] $\mathcal{D}=\cup_{i=1}^nE_i$ is a string such that $E_i$ is $\Ff_Y$-invariant and $Z(\Ff_Y,E_i)=2$ for any $i$. Moreover, in this case,
    \begin{enumerate}
        \item $X\ni x$ is either a smooth point or a cyclic quotient singularity and $\Ff$ is non-dicritical near $x$. In particular, $X\ni x$ is klt and $K_{\Ff}$ is $\Qq$-Cartier,
        \item $a(E_i,\Ff)=0$ for any $i$,
        \item $\pld(X\ni x,\Ff)=0$,  $B=0$, and $(X\ni x,\Ff)$ is canonical but not terminal, and
        \item $K_{\Ff}$ is Cartier near $x$.
    \end{enumerate}
    \item [(Case 4)] $\mathcal{D}(f)=\cup_{i=1}^3E_i\cup\cup_{j=1}^nF_j$ for some positive integer $n$ is the following:
\begin{center}
    \begin{tikzpicture}
                                                    
        % \draw (-1,0) ellipse (1.5 and 0.25);   %%\\short oval
        % \draw (0.75,0) circle (0.1);
        % \draw (-0.9,0) ellipse (2 and 0.45);   %% long oval
        % \draw (1.1,0)--(1.4,0);
        % \draw (1.5,0.6) circle (0.1);
        % \node [left] at (1.4,0.6) {\footnotesize$2$};
        % \draw (1.5,0.1)--(1.5,0.5);
        % \draw[fill=black] (1.5,0) circle (0.1);
        % \node [above] at (1.5,0.1) {\footnotesize$2$};
        % \draw (1.5,-0.6) circle (0.1);
        % \node [left] at (1.4,-0.6) {\footnotesize$2$};
        % \draw (1.5,-0.1)--(1.5,-0.5);
        % \draw (1.6,0)--(2,0);
         \draw (2.1,0) circle (0.1);
         \node [above] at (2.1,0.1) {\footnotesize$F_{n}$};
         \draw[dashed] (2.2,0)--(3.2,0);
        % \draw (2.7,0) circle (0.1);
        % \node [above] at (2.7,0.1) {\footnotesize$2$};
        % \draw (2.8,0)--(3.2,0);
         \draw (3.3,0) circle (0.1);
         \node [above] at (3.3,0.1) {\footnotesize$F_1$};
         \draw (3.4,0)--(3.8,0);

         \draw (3.9,0) circle (0.1);
         \draw (3.9,0.6) circle (0.1);
         \draw (3.9,-0.6) circle (0.1);
        % \draw (3.9,-1.2) circle (0.1);
         
         \draw (3.9,0.1)--(3.9,0.5);
         \draw (3.9,-0.1)--(3.9,-0.5);
        % \draw (3.9,-0.7)--(3.9,-1.1);
         
         \node [right] at (4,0) {\footnotesize$E_2$};
         \node [right] at (4,0.6) {\footnotesize$E_1$};
         \node [right] at (4,-0.6) {\footnotesize$E_3$};
        % \node [right] at (4,-1.2) {\footnotesize$2$};
    \end{tikzpicture},
\end{center}
such that
\begin{itemize}
    \item $\cup_{j=1}^nF_j$ is a Hirzebruch-Jung string and $Z(\Ff_Y,F_j)=2$ for any $j$, and
    \item $\mathcal{D}=\cup_{i=1}^3E_i$ is a Hirzebruch-Jung string such that $Z(\Ff_Y,E_1)=Z(\Ff_Y,E_3)=1,Z(\Ff_Y,E_2)=3$, $E_1^2=E_3^2=-2$, and $E_2^2\leq -2$.
\end{itemize}
Moreover, in this case,
\begin{enumerate}
    \item $X\ni x$ is a $D$-type singularity and $\Ff$ is non-dicritical near $x$. In particular, $X\ni x$ is klt and $K_{\Ff}$ is $\Qq$-Cartier near $x$,
    \item $a(E_1,\Ff)=a(E_3,\Ff)=\frac{1}{2}$, $a(E_2,\Ff)=0$, and $a(F_j,\Ff)=0$ for any $j$,  
    \item  $\pld(X\ni x,\Ff)=0$,  $B=0$,  and $(X\ni x,\Ff)$ is canonical but not terminal, and
    \item $2K_{\Ff}$ is Cartier near $x$.
\end{enumerate}
\item [(Case 5)] 
\begin{itemize}
    \item Either $\mathcal{D}=\cup_{i=1}^nE_i$ is a cycle such that each $E_i$ is $\Ff_Y$-invariant and $Z(\Ff_Y,E_i)=2$ for any $i$, or
    \item $\mathcal{D}=E_1$ is an $\Ff_Y$-invariant rational curve with a unique nodal singularity $x$, such that $x$ is a reduced singularity of $\Ff_Y$ and $Z(\Ff_Y,E_1)=0$.
\end{itemize} 
Moreover, in this case, 
\begin{enumerate}
    \item $X\ni x$ is an elliptic singularity, $\Ff$ is non-dicritical near $x$, and $K_{\Ff}$ is not $\Qq$-Cartier near $x$.
    \item $a(E_i,\Ff)=0$ for any $i$, and
    \item $\pld(X\ni x,\Ff)=0$, $B=0$, and $(X\ni x,\Ff)$ is num-canonical but not num-terminal.
\end{enumerate}

\item [(Case 6)] $\mathcal{D}=\cup_{i=1}^nE_i\cup D\cup_{j=1}^mF_j$ for some non-negative integers $m,n$ is the following:
\begin{center}
    \begin{tikzpicture}
                                                    
        % \draw (-1,0) ellipse (1.5 and 0.25);   %%\\short oval
        % \draw (0.75,0) circle (0.1);
        % \draw (-0.9,0) ellipse (2 and 0.45);   %% long oval
        % \draw (1.1,0)--(1.4,0);
        % \draw (1.5,0.6) circle (0.1);
        % \node [left] at (1.4,0.6) {\footnotesize$2$};
        % \draw (1.5,0.1)--(1.5,0.5);
        % \draw[fill=black] (1.5,0) circle (0.1);
        % \node [above] at (1.5,0.1) {\footnotesize$2$};
        % \draw (1.5,-0.6) circle (0.1);
        % \node [left] at (1.4,-0.6) {\footnotesize$2$};
        % \draw (1.5,-0.1)--(1.5,-0.5);
        % \draw (1.6,0)--(2,0);
         \draw (2.1,0) circle (0.1);
         \node [below] at (2.1,-0.1) {\footnotesize$E_{n}$};
         \draw[dashed] (2.2,0)--(3.2,0);
        % \draw (2.7,0) circle (0.1);
        % \node [above] at (2.7,0.1) {\footnotesize$2$};
        % \draw (2.8,0)--(3.2,0);
         \draw (3.3,0) circle (0.1);
         \node [below] at (3.3,-0.1) {\footnotesize$E_1$};
         \draw (3.4,0)--(3.8,0);

         %\draw (3.9,+1.5) circle (0.1);
         \draw (3.9,0) circle (0.1);
         %\draw (3.9,0.6) circle (0.1);
        % \draw (3.9,-0.6) circle (0.1);
        % \draw (3.9,-1.5) circle (0.1);
         \draw (4.5,0) circle (0.1);
         \node [below] at (4.7,-0.1) {\footnotesize$F_{1}$};
         \draw (5.7,0) circle (0.1);
         \node [below] at (5.7,-0.1) {\footnotesize$F_m$};

         %\draw [dashed](3.9,+0.7)--(3.9,+1.4);
         %\draw (3.9,0.1)--(3.9,0.5);
        % \draw (3.9,-0.1)--(3.9,-0.5);
        % \draw [dashed](3.9,-0.7)--(3.9,-1.4);
         \draw (4.0,0)--(4.4,0);
         \draw [dashed](4.6,0)--(5.6,0);

         %\node [right] at (4,1.5) {\footnotesize$F_{m_2}$};
         \node [below] at (3.9,-0.1) {\footnotesize$D$};
         %\node [right] at (4,0.6) {\footnotesize$F_{m_1+1}$};
        % \node [right] at (4,-0.6) {\footnotesize$F_{m_2+1}$};
        % \node [right] at (4,-1.5) {\footnotesize$F_{m}$};
    \end{tikzpicture},
\end{center}
such that
\begin{itemize}
    \item $D$ is not $\Ff_Y$-invariant and $\tang(\Ff_Y,D)=0$,
    \item either $n=0$ or $\cup_{i=1}^nE_i$ is an $\Ff_Y$-chain, 
    \item either $m=0$ or $\cup_{j=1}^mF_j$ is an $\Ff_Y$-chain, and
\end{itemize}
Moreover, in this case,
\begin{enumerate}
\item $\Ff$ is dicritical near $x$, and one of the following holds:
\begin{enumerate}
    \item [(Case 6.1)] $D$ is a rational curve. Then $X\ni x$ is a cyclic quotient singularity. In particular, $X\ni x$ is klt and $K_{\Ff}$ is $\Qq$-Cartier near $x$.
    \item [(Case 6.2)] $D$ is an elliptic curve and $m=n=0$.  Then $X\ni x$ is an elliptic singularity. In particular, $X\ni x$ is lc but not klt.
    \item [(Case 6.3)] $D$ is not a rational curve, and either $m>0$, or $n>0$, or $p_a(D)\geq 2$. Then $X\ni x$ is not lc,
\end{enumerate}
\item $a(D,\Ff)=-1$, $a(E_i,\Ff)=0$ for any $i$, and $a(F_j,\Ff)=0$ for any $j$,
\item $\pld(X\ni x,\Ff)=0$,  $B=0$, and $(X\ni x,\Ff)$ is num-lc but neither num-canonical nor num-klt, and
\item if $K_{\Ff}$ is $\Qq$-Cartier near $x$, then $K_{\Ff}$ is Cartier near $x$.
\end{enumerate}
\item[(Case 7)] $\mathcal{D}(f)=D\cup_{i=1}^n\cup_{j=1}^{r_i}E_{i,j}$ 
for some integer $n\geq 3$ and positive integers $r_1,\dots,r_n$, such that
\begin{itemize}
    \item $D$ is not $\Ff_Y$-invariant and $\tang(\Ff_Y,D)=0$, and
    \item for any $i$, $\cup_{j=1}^{r_i}E_{i,j}$ is an $\Ff_Y$-chain.
\end{itemize}
Moreover in this case,
\begin{enumerate}
    \item $\Ff$ is dicritical near $x$, and one of the following holds:
    \begin{enumerate}
        \item [(Case 7.1)] If $D$ is a rational curve, then $X\ni x$ is klt (resp. lc) if and only if $$
        \sum_{i=1}^n(1-\frac{1}{d_i})<2,~(\text{resp.} \leq 2),
        $$ where $d_i=\det(\mathcal{D}(\cup_{j=1}^{r_i}E_{i,j}))$ .
    \item [(Case 7.2)] If $D$ is not a rational curve, then $X\ni x$ is not klt,
    \end{enumerate}
    \item $a(D,\Ff)=-1$, and $a(E_{i,j},\Ff)=0$ for any $i,j$,
    \item $\pld(X\ni x,\Ff)=0$,  $B=0$, and $(X\ni x,\Ff)$ is num-lc but neither num-canonical nor num-klt, and
    \item  if $K_{\Ff}$ is $\Qq$-Cartier near $x$, then $K_{\Ff}$ is Cartier near $x$.
\end{enumerate}
\end{enumerate}
\end{thm}
\begin{proof}
By Lemma \ref{lem: add component worse sing}, $(X\ni x,\Ff)$ is num-lc. Now the main part of the theorem follows immediately from \cite[Theorem 2.4]{Che23}. More precisely, the only difference of the main part of our theorem from \cite[Theorem 2.4]{Che23} is that we do not assume that $K_{\Ff}$ is $\Qq$-Cartier near $x$. Nevertheless, since \cite[Proof of Theorem 2.4]{Che23} only relies on the structure of $\mathcal{D}$ and $\Ff_Y$, the same arguments of \cite[Proof of Theorem 2.4]{Che23} will provides the classification of $\mathcal{D}$ in our situation as well. We remark that there is a small difference for Case 6 comparing to \cite[Theroem 2.4]{Che23}: although \cite[Theorem 2.4(6)]{Che23} states that $\cup_{i=1}^nE_{n+1-i}$ is an $\Ff_Y$-chain in this case, it is actually $\cup_{i=1}^nE_i$ that is an $\Ff_Y$-chain. To see this, we may simply apply \cite[Theorem 2.4, (17) Claim]{Che23}.

In the following, we only prove the moreover part for each case of our theorem.

\medskip
\noindent(Case 1) Since all curves in $\mathcal{D}$ are $\Ff_Y$-invariant, $\Ff$ is non-dicritical near $x$. Since $\mathcal{D}$ is a chain of rational curves, $X\ni x$ is a cyclic quotient singularity, which implies (a). By Theorem \ref{thm: tang and Z formula}(2) and computing intersection numbers of $K_{\Ff_Y}+\sum_{i=1}^m(-a(E_i,\Ff))\cdot E_i$ with $E_i$, we get (b). (c) follows from (b). Since $\mathcal{D}$ is an $\Ff_Y$-chain, there exists a unique $\Ff_Y$-invariant curve $C_Y\not\subset\Supp\mathcal{D}$ on $Y$ which intersects $\mathcal{D}$, and $C_Y$ intersects $E_m$. We let $y:=C_Y\cap E_m$, then $\Ff_Y$ has a reduced singularity at $y$ and hence $C_Y+E_m$ is snc at $y$. Let $C=f_*C_Y$, then by \cite[Theorem 4.15(3)]{KM98} we know $(X\ni x, C)$ is plt. Therefore $C$ is normal at $x\in X$.

%and let $v$ be the vector field which generates $C_Y$ near $y$, then $f_*v$ is the vector field which generates $\Ff$ near $x$. Thus $f_*v$ does not vanish at $x$, so $x$ is not a singular point of $\Ff$.

\medskip

\noindent(Case 2) Since all curves in $\mathcal{D}$ are $\Ff_Y$-invariant, $\Ff$ is non-dicritical near $x$. Since $\mathcal{D}$ is a chain of rational curves, $X\ni x$ is a cyclic quotient singularity, which implies (a). By Theorem \ref{thm: tang and Z formula}(2) and computing intersection numbers of $K_{\Ff_Y}+\sum_{i=1}^3(-a(E_i,\Ff))\cdot E_i$ with $E_i$, we get (b). (c) follows from (b) and Lemma \ref{lem: add component worse sing}. (d) follows from (a) and (b).

\medskip
\noindent(Case 3) Since all curves in $\mathcal{D}$ are $\Ff_Y$-invariant, $\Ff$ is non-dicritical near $x$. Since $\mathcal{D}$ is a chain of rational curves, $X\ni x$ is a cyclic quotient singularity, which implies (a). By Theorem \ref{thm: tang and Z formula}(2), $K_{\Ff_Y}\cdot E_j$ for any $j$. Thus $(\sum_{i=1}^na(E_i,\Ff)\cdot E_i)\cdot E_j=0$ for any $j$. By the negativity lemma, we get (b). (c) follows from (b) and Lemma \ref{lem: add component worse sing}. (d) follows from (a) and (b).

\medskip

\noindent(Case 4) Since all curves in $\mathcal{D}$ are $\Ff_Y$-invariant, $\Ff$ is non-dicritical near $x$. Since all components of $\mathcal{D}$ is are rational curves, $X\ni x$ is a $D$-type singularity, which implies (a). Let $$G:=2\sum_{j=1}^na(F_j,\Ff)\cdot F_j+a(E_2,\Ff)\cdot(E_1+2E_2+E_3).$$ By Theorem \ref{thm: tang and Z formula}(2) and by computing intersection numbers, we know $a(E_1,\Ff)=a(E_3,\Ff)=\frac{a(E_2,\Ff)+1}{2}$,and $G\equiv_X 0$. By the negativity lemma, we get (b). (c) follows from (b) and Lemma \ref{lem: add component worse sing}. (d) follows from (a) and (b).

\medskip

\noindent(Case 5) Since all curves in $\mathcal{D}$ are $\Ff_Y$-invariant, $\Ff$ is non-dicritical near $x$. By classification of surface singularities, $X\ni x$ is an elliptic singularity. By \cite[Theorem IV.2.2]{McQ08}, $K_{\Ff}$ is not $\Qq$-Cartier near $x$. This implies (a). By Theorem \ref{thm: tang and Z formula}, $K_{\Ff_Y}\cdot E_i=0$ for any $i$. This implies (b). (c) follows from (b) and Lemma \ref{lem: add component worse sing}.

\medskip

\noindent(Case 6 (a-c)) Since $D$ is not $\Ff_Y$-invariant, $\Ff$ is dicritical near $x$. (a) follows from the classification of surface singularities. Let $$G:=\sum_{i=1}^na(E_i,\Ff)\cdot E_i+(1+a(D,\Ff))\cdot D+\sum_{j=1}^ma(F_j,\Ff)\cdot F_j.$$
By Theorem \ref{thm: tang and Z formula}, $G\equiv_X 0$. By the negativity lemma, we get (b). (c) follows from (b) and Lemma \ref{lem: add component worse sing}. We
%and , $(K_{\Ff_Y}+D)|_{\cup_{i=1}^nE_i\cup D\cup_{j=1}^mF_j}\sim 0$. By \cite[11.3.6 Lemma]{Kol+92}, we get (d).

\medskip

\noindent(Case 7 (a-c)) Since $D$ is not $\Ff_Y$-invariant, $\Ff$ is dicritical near $x$. (a) follows from the classification of surface singularities. Let $$G:=\sum_{i=1}^n\sum_{j=1}^{r_i}a(E_{i,j},\Ff)\cdot E_{i,j}+(1+a(D,\Ff))\cdot D.$$
By Theorem \ref{thm: tang and Z formula}, $G\equiv_X 0$. By the negativity lemma, we get (b). (c) follows from (b) and Lemma \ref{lem: add component worse sing}.  By Lemma \ref{lem: restrction tang}, we have $(K_{\Ff_Y}+D)|_{D}=0$. 

\medskip

\noindent{(Case 6(d) and Case 7(d))} By Lemma \ref{lem: restrction tang}, $(K_{\Ff_Y}+D)|_{D}\sim0$. We let $C$ be the reduced $f$-exceptional divisor. Then since $\mathcal{D}$ is a tree, $E_i,F_j$ are smooth rational curves in Case 6, and $E_{i,j}$ are smooth rational curves in Case 7, by \cite[Proposition 7.5.4]{Liu06},
$$
H^1(C,\Oo_C)=p_a(C)=p_a(D)=H^1(D,\Oo_D).
$$
Therefore, \cite[Theorem 7.5.19]{Liu06} implies that the canonical homomorphism
$$
\Pic^0(C)\to \Pic^0(D)
$$
is an isomorphism. Since $K_{\Ff_Y}+D\equiv_X0$ and $(K_{\Ff_Y}+D)|_{D}\sim0$, we have $(K_{\Ff_Y}+D)|_C\sim0$. 

If $K_\Ff$ is $\Qq$-Cartier, then $f^*K_\Ff=K_{\Ff_Y}+D$. 
Next we can choose a contractible stein neighborhood $V$ of $x\in X$ such that $U=f^{-1}(V)$ deformation retracts to $C$, then $\Pic(V)$ is trivial by the following exact sequence
$$
\cdots\to H^1(V,\Oo_V)\to H^1(V,\Oo_V^*)\simeq\Pic(V)\to H^2(V,\Zz)\to\cdots.
$$
Moreover, the canonical homomorphisms $H^i(U,\ZZ)\to H^i(C,\Zz)$ induced by the inclusion are isomorphisms. Suppose that $rK_\Ff$ is Cartier, then $rK_\Ff|_V$ is trivial and hence $r(K_{\Ff_Y}+D)|_U$ is also trivial. This implies that $(K_{\Ff_Y}+D)|_U$ is a torsion in $\Pic(U)$, then by \cite[11.3.6 Lemma]{Kol+92} we know that $(K_{\Ff_Y}+D)|_U\sim 0$ and hence its pushforward $K_\Ff$ is also trivial in $\Pic(V)$. Therefore $K_\Ff$ is Cartier, and we get (d) for both cases.
\end{proof}

\begin{cor}\label{cor: dlt foliated surface}
Let $(X,\Ff,B)$ be a dlt (cf. \cite[Definition 3.6]{CS21}) foliated triple such that $\dim X=2$ and $\rk\Ff=1$. Then for any closed point $x\in X$,
\begin{enumerate}
    \item if $x$ is a singular point of $X$, then $(X\ni x,\Ff,B)$ is as in Case 1 of Theorem \ref{thm: classification of foliation intro}. In particular, $x$ is a non-singular point of $\Ff$, $\mld(X\ni x,\Ff,B)>0$, and $x$ is a cyclic quotient singularity of $X$, and
    \item if $x$ is a non-singular point of $X$, then one of the following cases hold:
    \begin{enumerate}
        \item $x$ is a non-singular point of $\Ff$.
        \item $x$ is a reduced singularity of $\Ff$.
        %\item $(X\ni x,\Ff,B)$ is as in Case 3 of Theorem \ref{thm: classification of foliation intro} such that $\det(\mathcal{D})=1$. In particular, $x\not\in\Supp B$ and $\mld(X\ni x,\Ff)=0$.
    \end{enumerate}
\end{enumerate}
\end{cor}
\begin{proof}
It immediately follows from Theorem \ref{thm: classification of foliation intro}.
\end{proof}

\subsection{Examples}

\begin{exrem}\label{exrem: counterexample cs}
We remark that for each singularity listed in Theorem \ref{thm: classification of foliation intro}, there are corresponding examples. Indeed, it is very easy to construct those example by considering the foliation induced by the natural $\mathbb P^1$-bundle structure $\mathbb P(E)\rightarrow E$ for some curve $E$. Examples for Case 1-4 of Theorem \ref{thm: classification of foliation intro} can be constructed by constructing a sequence of blow-ups along a general fiber of $\mathbb P(E)\rightarrow E$. Examples for Case 6 and 7 of Theorem \ref{thm: classification of foliation intro} can be constructed by taking blow-ups along a negative section of $\mathbb P(E)\rightarrow E$ or $\mathbb{P}^1\times\mathbb{P}^1$, then keep blowing up at the reduced singularities and blow-down the strict transform of the negative section along with some chain of rational curves. Examples for Case 5 can be constructed by considering a family of elliptic curves $X\rightarrow Z$ with a singular fiber $X_0$, blowing up a point on $X_0$, and contract the strict transform of $X_0$. We also remark that Examples for Cases 1-4 can be found in \cite[Fact I.2.4]{McQ08}.

In particular, the examples for Case 7.1 provide a negative answer to a question of P. Cascini and C. Spicer \cite[Question 3]{CS} on whether rational lc foliated surface germs are quotient singularities: when $\sum_{i=1}^n(1-\frac{1}{d_i})\geq 2$, $X\ni x$ is not klt, so $X\ni x$ is no longer a quotient singularity. All we need to guarantee is that $x\in X$ is a rational singularity simultaneously. Fortunately, examples arise from \cite[Chapter 3, Exercise 13]{Fri98}, where we have $\mathcal{D}(f)=D\cup_{i=1}^{4}E_i$ and $x\in X$ is a rational singularity as long as $D^2\le -3$.

%However, $\mathcal{D}(f)$ is a tree of rational curves, so $x$ is a rational singularity. 
\end{exrem}

Next we present a example showing that Grauert-Riemenschneider type vanishing theorem fails for surface foliations in general. This example may be well-known to experts, but can also be a good exercise for beginners in the foliation theories.

\begin{ex}\label{ex: case 67 not qcartier}
 Let $C$ be a smooth curve of genus $g\ge2$ and $S:=C\times C$. Let $p_i:S\to C~,i=1,2$ be the corresponding projections. If $\Delta:E\to S$ is the diagonal morphism, then $E$ is isomorphic to $C$ and we have $E^2=2-2g<0$. Let $\Ff$ be the foliation on $S$ determines by 
 $$
 0\to T_{\Ff}:=p^*_1T_C\to T_S\to p^*_2T_C\to 0,
 $$
 which is exactly the foliation induced by the fibration $p_2:S\to C$ and $K_\Ff=p^*_1K_C$ . By \cite[3.0 Theorem]{Kee99}, we have
 \begin{itemize}
     \item $p^*_iK_C+E$ is nef and big but not semi-ample, and
     \item $K_S+2E$ is semi-ample and defines a birational morphism $f:S\to Z$ which only contracts $E$.
 \end{itemize}
Let $\Ff_Z:=f_*\Ff$ be the pushforward foliation on $Z$ which is determined by $\Ff_Z|_{Z\backslash\{z\}}=\Ff|_{X\backslash E}$, where $z=f(E)$. Then $K_{\Ff_Z}=f_*K_\Ff$ and it is easy to see that $(Z,\Ff_Z)$ is a num-lc foliated surface (Case 6 or 7 of Theorem \ref{thm: classification of foliation intro}). Notice that $\Ff_Z$ is dicritical at $z\in Z$ and the minimal resolution of $\Ff_Z$ is $f:S\to Z$. 
 %Let $\Ff$ be the foliation on $S$ determines by 
 %$$
 %0\to T_{\Ff}:=p^*_1\Oo_C(K_C)\to T_S\to p^*_2\Oo_C(K_C)\to 0,
 %$$
 %which is exactly the foliation induced by the fibration $p_2:S\to C$. Let $\Ff_Z:=f_*\Ff$ be the pushforward foliation on $Z$ which is determined by $\Ff_Z|_{Z\backslash\{z\}}=\Ff|_{X\backslash E}$, where $z=f(E)$. Then $K_{\Ff}=L_1$ and $K_{\Ff_Z}=f_*K_\Ff$. It is easy to see that $(Z,\Ff_Z)$ is a num-lc foliated surface (Case 6 or 7 of Theorem \ref{thm: classification of foliation intro}), $\Ff_Z$ is dicritical at $z\in Z$, and the minimal resolution of $\Ff_Z$ is $f:S\to Z$.

 We claim that 
 \begin{enumerate}
     \item $(K_{\Ff}+E)|_E\sim 0$ but $K_{\Ff_Z}=f_*K_\Ff$ is not $\Qq$-Cartier. 
     \item $K_\Ff+E$ is not semi-ample over $Z$ and $R^1f_*\Oo_S(K_\Ff)\neq 0$. In particular, Grauert-Riemenschneider type vanishing theorem fails for $f: (S,\Ff)\rightarrow (Z,\Ff_Z)$.
 \end{enumerate}
 For (1), notice that $(p^*_iK_C)|_E=(p_i\circ\Delta)^* K_C=K_E$. By adjunction we have 
 $$
 K_E=(K_S+E)|_E=(p_1^*K_C+p_2^*K_C+E)|_E=2K_E+E|_E,
 $$
 therefore $\Oo_S(E)|_E\sim -K_E$ and $(K_\Ff+E)|_E=(p_1^*K_C+E)|_E\sim 0$. 
 
 If $K_{\Ff_Z}$ is $\Qq$-Cartier, then $f^*K_{\Ff_Z}=K_\Ff+aE$ and $a$ must be 1 by the previous statement. Let $C'$ be any irreducible curve on $Z$ and let $\bar C$ be its strict transform on $S$. Since $\bar C\neq E$, $\bar C\cdot E>0$. Then we can see that $K_{\Ff_Z}\cdot C'=(K_\Ff+E)\cdot \bar C>0$. Indeed, either
 \begin{itemize}
     \item $p_1(\bar C)=C$ so that $K_\Ff\cdot\bar C>0$, or
     \item $p_1(\bar C)$ is a single point so that $\bar C\cdot E>0$.     
 \end{itemize}
Therefore the big divisor $K_{\Ff_Z}$ is actually ample, which implies that $K_{\Ff}+E=f^*K_{\Ff_Z}$ is semi-ample and we reach a contradiction. 

Next we prove (2). If $K_\Ff+E$ is semi-ample over $Z$, then there exists a morphism $g:S\to Y$ over $Z$ defined by $K_\Ff+E$. Since $K_\Ff+E$ is not ample over $Z$, $g$ is not an isomorphism so that $Y=Z$. Therefore $K_\Ff+E$ is a pullback of a $\Qq$-Cartier divisor on $Z$ and this divisor is necessarily the pushforward $f_*(K_\Ff+E)=K_{\Ff_Z}$, which contradicts (1). Hence $K_\Ff+E$ is not semi-ample over $Z$. %, in particular, is not base point free over $Z$. Considering the long exact sequence %induced by
Consider the long exact sequence
$$
0\to f_*\Oo_S(K_\Ff)\stackrel{\alpha}{\longrightarrow}f_*\Oo_S(K_\Ff+E)\stackrel{\beta}{\longrightarrow}f_*\Oo_E\stackrel{\gamma}{\longrightarrow}R^1f_*\Oo_S(K_\Ff)\to\cdots.
$$
Since $f_*\Oo_E=H^0(E,\Oo_E)\simeq\Cc$, we only need to show that $\alpha$ is surjective, so that $\beta$ is zero, $\gamma$ is injective, and $R^1f_*\Oo_S(K_\Ff)\neq 0$. 

If $\alpha$ is not surjective, then there exists an effective divisor $D\sim_{Z}K_\Ff+E$ such that $E$ is not in the support of $D$. Notice that $D|_E\sim 0$ so $D$ must be disjoint from $E$. This is impossible since then $f_*D$ is a Cartier divisor on $Z$, so $f_*(K_\Ff+E)\sim_Z D$ is also a Cartier divisor on $Z$, contradicting (1). 
\end{ex}

\begin{rem}
We note that the Grauert-Riemenschneider type vanishing theorem for foliated surfaces has been established for canonical singularities according to \cite[Theorem 5]{HL21}. Furthermore, it has been extended to ``good log canonical" singularities as defined in \cite[Definition 5.1]{Che23} through \cite[Theorem 0.3]{Che23}.

A. Langer has informed us that M. Lupinski \cite{Lup21,Lup} has independently found another counterexample to the Grauert-Riemenschneider type vanishing theorem for contractions between foliated surfaces with lc singularities by considering the minimal resolution of the foliation $\Ff$ as described in \cite[Example I.2.5]{McQ08}. Specifically, in \cite{Lup21}, Lupinski has discovered the minimal resolution $f: (Y,\Ff_Y) \rightarrow (X,\Ff)$ of $\Ff$, while \cite[1.2 Grauert-Riemenschneider type vanishing theorem]{Lup} demonstrates that the Grauert-Riemenschneider type vanishing theorem fails for the morphism $f: (Y,\Ff_Y) \rightarrow (X,\Ff)$, i.e. $R^1f_*\mathcal{O}_Y(K_{\Ff_Y}) \neq 0$.

However, it remains an open question whether the Grauert-Riemenschneider type vanishing theorem holds for foliations with canonical singularities in higher dimensions (see \cite[Question 6]{HL21}).
\end{rem}

%so that in particular $f_*\Oo_S(K_\Ff+E)$ cannot generate $\Oo_S(K_\Ff+E)$. 

\section{Proof of the main theorems}\label{sec: proof of main theorems}

The following theorem is important for the proof of our main results.

\begin{thm}\label{thm: surface from foliation to plt pair}
Let $(X\ni x,\Ff)$ be a surface foliated germ such that either both $X$ and $\Ff$ are smooth near $x$ or $(X\ni x,\Ff)$ is as in Case 1 of Theorem \ref{thm: classification of foliation intro}. Then:
\begin{enumerate}
    \item There exists a unique $\Ff$-invariant irreducible curve $L$ passing through $x$.
    \item For any $B\geq 0$ on $X$ and any prime divisor $E$ over $X\ni x$,
    $$a(E,\Ff,B)=a(E,X,B+L)+1.$$
\end{enumerate}
Note that although $L$ may not be algebraic, it is at least locally analytically well-defined.
\end{thm}
\begin{proof}
Let $f: Y\rightarrow X$ be the minimal resolution of $\Ff\ni x$, $\Ff_Y:=f^{-1}\Ff$, and $E_1,\dots,E_m$ $f$-exceptional prime divisors, such that either $m=0$ or $\cup_{i=1}^mE_i$ is an $\Ff_Y$-chain. 

(1) If $X$ and $\Ff$ are both smooth near $x$ then there is nothing left to prove. So we may assume that $(X\ni x,\Ff)$ is as in Case 1 of Theorem \ref{thm: classification of foliation intro}. Then it follows from Theorem \ref{thm: classification of foliation intro} (Case 1.d) and we let $L$ be that curve.  %there exists a closed point $y\in E_m$ such that $y\not=E_{m-1}\cap E_m$ and $Z(f^{-1}\Ff,E_m)=1$. Since $\Ff_Y$ has at most reduced singularities, there exists a unique $\Ff_Y$-invariant irreducible curve $L_Y$ on $Y$ such that $y=L_Y\cap E_m$. We may let $L:=f_*L_Y$.

(2) We only need to prove the case when $B=0$. For any prime divisor $E$ over $X\ni x$, there exists a sequence of blow-ups 
$$Y_n\xrightarrow{h_n} Y_{n-1}\xrightarrow{h_{n-1}}\cdots\xrightarrow{h_2} Y_1\xrightarrow{h_1} Y_0:=Y,$$
such that
\begin{itemize}
    \item $E$ is on $Y_n$ but not on $Y_{n-1}$,
    \item $\Ff_0:=\Ff_Y$ and $\Ff_i:=h_i^{-1}\Ff_{i-1}$ for each $i$,
    \item $F_i:=\Exc(h_i)$ is a prime $h_i$-exceptional divisor for each $i$, and
    \item $h_i$ is the blow-up of a closed point $y_{i-1}\in Y_{i-1}$ such that $y_{i-1}$ is contained in the union of the strict transforms of $F_1,\dots,F_{i-1},E_1,\dots,E_m$ and $L$.
\end{itemize}
We prove (2) by applying induction on the number $n$. When $n=0$, (2) follows from Theorem \ref{thm: classification of foliation intro}(Case 1.c). Suppose that $n>0$. There are two cases.

\medskip

\noindent\textbf{Case 1}. $y_{n-1}$ is contained in exactly two curves $C_1,C_2$ of the strict transforms of $F_1,\dots,F_{i-1},E_1,\dots,E_m$ and $L$. By the induction, $a(C_1,\Ff,0)=a(C_1,X,L)+1$ and $a(C_2,\Ff,0)=a(C_2,X,L)+1$. So
\begin{align*}
  a(E,\Ff,0)&=a(E,\Ff_{n-1},-a(C_1,\Ff,0)C_1-a(C_2,\Ff,0)C_2)=a(C_1,\Ff,0)+a(C_2,\Ff,0)\\
  &=a(C_1,X,L)+a(C_2,X,L)+2\\
  &=a(C_1,X,L)\mult_EC_1+a(C_2,X,L)\mult_EC_2+a(E,Y_{n-1},0)+1\\
  &=a(E,Y_{n-1},-a(C_1,X,L)C_1-a(C_2,X,L)C_2)+1=a(E,X,L)+1
\end{align*}
and we are done.

\medskip

\noindent\textbf{Case 2}. $y_{n-1}$ is contained in exactly one curve $C$ of the strict transforms of $F_1,\dots,F_{i-1},E_1,\dots,E_m$ and $L$. By the induction, $a(C,\Ff,0)=a(C,X,L)+1$. So
\begin{align*}
  a(E,\Ff,0)&=a(E,\Ff_{n-1},-a(C,\Ff,0)C)=a(C,\Ff,0)+1\\
  &=a(C,X,L)+2=a(E,Y_{n-1},-a(C,X,L)C)+1=a(E,X,L)+1
\end{align*}
and we are done.
\end{proof}

\subsection{Boundedness of complements}

The following lemma might be well-known to experts, and is a standard technique when one are tackling between the algebraic setting and the analytic setting. 

\begin{lem}\label{lem: analytic N-complement to algebraic N-complement}
    Let $(X\ni x,B)$ be an lc germ with $\Qq$-factorial singularity. Assume that there exists an N-complement $B^+_{an}$ for $(X_{an}\ni x,B_{an})$ in the analytic setting, then there exists an N-complement $B^+$ for $(X\ni x,B)$ in the algebraic setting.
\end{lem}
\begin{proof}
We may assume that $x\in X\sim \Spec (R,\mathfrak{m})$ and $x\in X_{an}\sim\Spec(\hat{R},\hat{\mathfrak{m}})$, where $R$ is a localization of certain finite generated $\Cc$-algebra. Let 
$$
D:=N\lf B\rf+\lf(N+1)\{B\}\rf,~D_{an}:=N\lf B_{an}\rf+\lf(N+1)\{B_{an}\}\rf,
$$
then $\Oo_X(-NK_X-D)$, $\Oo_{X_{an}}(-NK_{X_{an}}-D_{an})$ are  rank one reflexive sheaves on $X$, $X_{an}$. Moreover, we have the relation:
$$
\Oo_{X_{an}}(-NK_{X_{an}}-D_{an})=\Oo_X(-NK_X-D)\otimes_R\hat{R}.
$$
Notice that $\Oo_X(-N_X-D)\cong\{f\in K(R)~|\Div(f)-NK_X-D\ge0\}$ and 
$$
\Oo_{X_{an}}(-NK_{X_{an}}-D_{an})\cong\{f\in K(\hat R)~|\Div(f)-NK_{X_{an}}-D_{an}\ge0\}
$$
Therefore by our assumption the effective Weil divisor $NB^+_{an}-D_{an}$ corresponds to certain non-zero element $f$ in $\Oo_{X_{an}}(-NK_{X_{an}}-D_{an})$ such that 
$$
\Div(f)-NK_{X_{an}}-D_{an}=NB^+_{an}-D_{an}.
$$
In particular, $\Oo_X(-NK_X-D)$ is not zero since $R\to\hat R$ is faithfully flat. Recall that
$$
|-NK_X-D|=\{\Div(f)-NK_X-D\ge0~|f\in K(R)\}.
$$
Let $\phi: X'\to X$ be a log resolution such that 
\begin{enumerate}
    \item $\phi^*|-NK_X-D|=F+|M|$, where $F$ is the fixed part and $M$ is base point free.
    \item $\phi_*^{-1}B+\Exc(\phi)+F$ is snc.
\end{enumerate}
Here we explain why the pullback in (1) makes sense: Since $-NK_X-D$ is $\Qq$-Cartier, we can define $D':=\phi^*(-NK_X-D)$ as a $\Qq$-divisor on $X'$, then we have
$$
\Div(f)-NK_X-D\ge0\Leftrightarrow \Div(f)+D'\ge0
$$
and we can set $|D'|:=\{\Div(f)+D'\ge0~|f\in K(R)\}$. Then we can define the fixed part $F$ of $|D'|$ as an effective $\Qq$-divisor, possibly replaced by a higher resolution, we can assume that $|D'|-F$ is base point free.

Next we prove that the fixed part of $|D'_{an}|$ is $F_{an}$. It is easy to see that $\mathrm{Fix}|D'_{an}|\le F_{an}$, so we just need to check $\mult_{\tilde{F}_i}|D'|=\mult_{F_i}|D'_{an}|$ for any irreducible component $F_i$ of $F$ and any irreducible component $\tilde{F}_i$ of $(F_i)_{an}$. Notice that $F_i$ corresponds to a discrete valuation $\nu_i$ of $K(R)$ and $\tilde{F}_i$ corresponds to a discrete valuation $\tilde{\nu}_i$ of $K(\hat R)$ that extends $\nu_i$. Then we have
$$
\mult_{F_i}|D'|=\min\{\mult_{F_i}(\Div(f)+D')~|~f\in\Oo_X(-NK_X-D)\subset K(R)\}.
$$
Since $\Oo_X(-NK_X-D)$ is coherent, we can find finite generators $f_j$ as an $R$-module, which will also be generators of $\Oo_{X_{an}}(-NK_{X_{an}}-D_{an})=\Oo_X(-NK_X-D)\otimes_R\hat R$ as an $\hat R$-module. By the properties of valuations, we get
\begin{align*}
\mult_{\tilde F_i}|D'_{an}|&=\min_j\{\mult_{\tilde F_i}(\Div(f_j)+D'_{an})\}\\
                    &=\min_j\{\tilde\nu_i(f_j)\}-\tilde\nu_i(D'_{an})\\
                    &=\min_j\{\nu_i(f_j)\}-\nu_i(D')\\
                    &=\mult_{F_i}|D'|
\end{align*}
Notice that $\phi_{an}:X'_{an}\to X_{an}$ is also a log resolution. Since $|M|,|M_an|$ is base point free, a general element in $|\Oo_X(-NK_X-D)|$ (resp. $|\Oo_{X_{an}}(-NK_{X_{an}}-D_{an})|$) induces an N-complement iff some inequality $\mult_{F_i}F\le a_i$ (resp. $\mult_{\tilde F_i}F_{an}\le a_i$) hold for all $i$. Therefore we are done by the existence of $B^+_{an}$.

\end{proof}

\begin{rem}\label{rem: more N-complement from An to Alg}
The general principle is that an algebraic linear system $|D'|$ behaves as good as its corresponding analytic linear system $|D'_{an}|$ over the germ $X\ni x$. This attributes to the relation between $\Oo_X(D')$ and $\Oo_{X_{an}}(D'_{an})$ at $x\in X$. One can easily check that the proof works for ``($\epsilon,$ N)-complement" since the only difference would be the inequalities conditions for the coefficients of $F$ and $F_{an}$.

\end{rem}

\begin{proof}[Proof of Theorem \ref{thm: lc complement foliated surface local}]
If $\pld(X\ni x,\Ff)=0$, then $\epsilon=0$. By Theorem \ref{thm: classification of foliation intro}, we may take $N=2$ and we are done. So we may assume that $\pld(X\ni x,\Ff)>0$.  By Theorem \ref{thm: classification of foliation intro}, $(X\ni x,\Ff,B)$ is either smooth or as in Case 1 of Theorem \ref{thm: classification of foliation intro}. In particular, $X\ni x$ is either smooth or a cyclic quotient singularity. By Theorem \ref{thm: surface from foliation to plt pair}, there exists a unique prime $\Ff$-invariant curve $L$ passing through $x$, such that $(X\ni x,B+L)$ is (analytically locally) $\epsilon$-lc.

Since $X\ni x$ is klt, by \cite[Theorem 1.1]{CH21}, there exists a positive integer $N$ depending only on $\Ii$, such that analytically locally, $(X\ni x,B+L)$ has an $(\epsilon,N)$-complement $(X\ni x,\tilde B^++L)$, and if $\bar\Ii\subset\mathbb Q$, then we may take $\tilde B^+\geq B$. Here we remark that since $L$ may not be algebraic, $\tilde B^+$ may also not be algebraic. We also remark that although \cite[Theorem 1.1]{CH21} only deals with algebraic pairs, the same lines of the proof works in the analytic setting.

By Theorem \ref{thm: surface from foliation to plt pair}(2), $(X\ni x,\Ff,\tilde B^+)$ is analytically locally $\epsilon$-lc. Let $f: Y\rightarrow X$ be the minimal resolution of $X\ni x$ and $E$ the reduced $f$-exceptional divisor. Let $L_Y:=f^{-1}_*L$ and $K_{\Ff_Y}+\tilde B^+_Y:=f^*(K_{\Ff}+\tilde B^+)$. By Theorem \ref{thm: surface from foliation to plt pair}(2),
$$K_Y+\tilde B_Y^++L_Y=f^*(K_X+\tilde B^++L)-E.$$
Thus $N(K_Y+\tilde B_Y^++L_Y)$ is Cartier over a neighborhood of $x$, so $N(K_{\Ff_Y}+\tilde B_Y^+)$ is Cartier over a neighborhood of $x$. Since $X\ni x$ is a cyclic quotient singularity, $N(K_{\Ff}+\tilde B^+)$ is Cartier near $x$. Thus, analytically locally, $(X\ni x,\Ff,\tilde B^+)$ is an $(\epsilon,N)$-complement of $(X\ni x,\Ff,B)$. 

Notice that for any resolution (of singularity) $\phi:X'\to X\in x$, assume that $(X\ni x,\Ff, C)$ is $\epsilon$-lc for some boundary $C$, then $(X',\Ff_{X'},C'+A)$ is also sub-$\epsilon$-lc over $x\in X$, where $K_{\Ff_{X'}}+C'=\phi^*(K_\Ff+C)$ and $A$ is a general element of a base point free system $|M|$ over $X\ni x$. This is because the exceptional curves are all invariant, $A\cap \phi^{-1}_*C=\emptyset$ over (a neighborhood of) $x$, and $A$ intersect the exceptional locus transversally at general smooth points of $\Ff_{X'}$. Since $K_{\Ff}$ is algebraic, we can apply exactly the same proof of Lemma \ref{lem: analytic N-complement to algebraic N-complement} by simply replacing $K_X$ with $K_{\Ff}$. As a consequence, we prove the existence of algebraic $(\epsilon,N)$-complement of $(X\ni x,\Ff,B)$.

%Since $X\ni x$ is either smooth or a cyclic quotient singularity, we have $(X\ni x)\cong\mathbb C^2/G$ for some cyclic group $G$, and we may identify $\tilde B^+$ with $(\frac{1}{N}(s=0))/G$ where $s$ is a formal power series in $\mathbb C[[x_1,x_2]]$. Let $l$ be a sufficiently large positive integer, $s_l$ the $l$-th truncation of $s$ (cf. \cite[Definition B.5]{HLL22}), and let $B^+:=(\frac{1}{N}(s_l=0))/G$. Then $B^+$ is an algebraic $\Qq$-divisor and $(X\ni x,\Ff,B^+)$ is an (algebraic) $(\epsilon,N)$-complement of $(X\ni x,\Ff,B)$.
\end{proof}

\begin{thm}\label{thm: 1 2 complement foliation nonintro}
Let $(X\ni x,\Ff)$ be a foliated lc surface germ such that $\rk\Ff=1$. Then $(X\ni x,\Ff)$ has a $2$-complement. Moreover, $(X\ni x,\Ff)$ has a $1$-complement if and only if $(X\ni x,\Ff)$ is not of Case 2 or Case 4 of Theorem \ref{thm: classification of foliation intro}.
\end{thm}

\begin{proof}[Proof of Theorem \ref{thm: 1 2 complement foliation nonintro}]
It suffices to consider the case $x\in X$ is not smooth, otherwise we are done.  By Theorem \ref{thm: classification of foliation intro}, if $\pld(X\ni x,\Ff)=0$ then $(X\ni x,\Ff)$ is a $2$-complement of itself, and is a $1$-complement of itself if and only if $(X\ni x,\Ff)$ is not of Case 2 or Case 4 of Theorem \ref{thm: classification of foliation intro}. So we may assume that $\pld(X\ni x,\Ff)>0$ and %If $X$ and $\Ff$ are non-singular near $x$, then $(X\ni x,\Ff)$ is a $1$-complement of itself, so we may assume that
$(X\ni x,\Ff)$ is of Case 1 of Theorem \ref{thm: classification of foliation intro}. Let $f: Y\rightarrow X$ be the minimal resolution of $X\ni x$ and let $E_1$ be the unique $f^{-1}\Ff$-invariant $f$-exceptional curve such that $Z(f^{-1}\Ff,E_1)=1$. We let $C_Y$ be any non-singular non-$f^{-1}\Ff$-invariant curve such that $C_Y$ intersects $E_1$ and $C_Y\cup\Exc(f)$ is snc. By computing intersection numbers, we know that $(X\ni x,\Ff,C:=f_*C_Y)$ is a $1$-complement of $(X\ni x,\Ff)$.
\end{proof}

\begin{proof}[Proof of Theorem \ref{thm: 1 2 complement foliation}]
It follows from Theorem \ref{thm: 1 2 complement foliation nonintro}.
\end{proof}

\subsection{The local index theorem}

\begin{proof}[Proof of Theorem \ref{thm: surface index theorem local}]
By Theorem \ref{thm: 1 2 complement foliation}, there exists a positive integer $I$ depending only on $a$ and $\Ii$ such that $(X\ni x,\Ff,B)$ has a monotonic $(a,I)$-complement $(X\ni x,\Ff,B^+)$. Then $B=B^+$, so $I(K_{\Ff}+B)$ is Cartier near $x$.
\end{proof}

\subsection{Set of mlds}

\begin{proof}[Proof of Theorem \ref{thm: set of foliated surface mlds}]
First we show ``$\subset$". Let $(X\ni x,\Ff,B)$ be a foliated triple such that $\dim X=2,~\rk\Ff=1$, and $B\in\Ii$. If $\pld(X\ni x,\Ff)=0$, then $\mld(X\ni x,\Ff,B)=0$ and we are done. So we may assume that  $\pld(X\ni x,\Ff)>0$. By Theorem \ref{thm: classification of foliation intro} (Case 1.d), there exists a unique prime $\Ff$-invariant curve $L$ passing through $x$, such that $\mld(X\ni x,B+L)=\mld(X\ni x,\Ff,B)$. We let $K_L+B_L:=(K_X+B+L)|_L$. Then $$B_L\in\left\{\frac{n-1+\sum c_i\gamma_i}{n}\mid n\in\mathbb N^+,c_i\in\mathbb N,\gamma_i\in\Ii\right\}\cap [0,1],$$
so $\mld(L\ni x,B_L)\in \{0,\frac{1-\sum c_i\gamma_i}{n}\mid n\in\mathbb N^+,c_i\in\mathbb N,\gamma_i\in\Ii\}\cap [0,1].$
By precise inversion of adjunction for surfaces,
$$\mld(L\ni x,B_L)=\mld(X\ni x,B+L)=\mld(X\ni x,\Ff,B),$$
so  $\mld(X\ni x,\Ff,B)\in\left\{0,\frac{1-\sum c_i\gamma_i}{n}\mid n\in\mathbb N^+,c_i\in\mathbb N,\gamma_i\in\Ii\right\}\cap [0,1].$ We remark that although $L$ may not be algebraic, we may still apply adjunction and inversion of adjunction to $L$ (cf. \cite[16.6 Proposition]{Kol+92}).

Now we show ``$\supset$". By Theorem \ref{thm: classification of foliation intro}, $0=\mld(X\ni x,\Ff)$ for some $(X\ni x,\Ff)$ such that $\dim X=2$ and $\rk\Ff=1$. Let $\Ff_0$ be the foliation induced by the natural fibration structure of $X_0:=\mathbb P^1\times\mathbb P^1\rightarrow Z:=\mathbb P^1$, $x_0\in X_0$ a closed point, $F$ the fiber of $X_0\rightarrow Z$ containing $x_0$, and $B_{i,j,0}$ general horizontal$/Z$ smooth rational curves. We blow-up the intersection of (the birational transform of) $F$ with the inverse image of $x_0$ $n$ times and get a contraction $h_n: X_n'\rightarrow X_0$. We let $F_n:=(h_n^{-1})_*F_0$, $B_{i,j,n}':=(h_n^{-1})_*B_{i,j,0}$, and $\Ff_n':=h_n^{-1}\Ff_0$. We let $g_n: X_n'\rightarrow X_n$ be the contraction of $F_n$, $x_n:=\Center_{X_n}F_n'$, $B_{i,j,n}:=(g_n)_*B_{i,j,n}'$, and $\Ff_n:=(g_n)_*\Ff_n'$. Let $B:=\sum_i\sum_{i=1}^{c_i}\gamma_iB_{i,j,n}$. Then 
$$a(F_n,\Ff_n,B)=\mld(X_n\ni x_n,\Ff_n,B)=\frac{1-\sum c_i\gamma_i}{n}$$
if $\sum c_i\gamma_i\leq 1$, and $a(F_n,\Ff_n,B)=\mld(X_n\ni x_n,\Ff_n,B)=-\infty$ otherwise. Thus ``$\supset$" holds.
\end{proof}

\begin{proof}[Proof of Corollaries \ref{cor: set of foliated surface mlds no boundary} and \ref{cor: acc mld foliation surface reprove}]
They are immediately implied by Theorem \ref{thm: set of foliated surface mlds}.
\end{proof}

\begin{proof}[Proof of Theorem \ref{thm: uniform bound of mlds}]
If $\pld(X\ni x,\Ff)=0$, then we may take $l=0$ and we are done. So we may assume that $\pld(X\ni x,\Ff)>0$. By Theorem \ref{thm: classification of foliation intro}, $(X\ni x,\Ff,B)$ is either smooth or as in Case 1 of Theorem \ref{thm: classification of foliation intro}. By Theorem \ref{thm: classification of foliation intro}(Case 1.d), there exists a unique prime $\Ff$-invariant curve $L$ passing through $x$, such that $(X\ni x,B+L)$ is lc. By \cite[Theorem 1.2]{HL22b}, there exists a positive integer $l$ depending only on $\Ii$ and a prime divisor $E$ over $X\ni x$, such that $a(E,X,B+L)=\mld(X\ni x,B+L)$ and $a(E,X,0)\leq l$. By Theorem \ref{thm: surface from foliation to plt pair}(2), $a(E,\Ff,B)=\mld(X\ni x,\Ff,B)$ and $a(E,\Ff,0)=a(E,X,L)\le a(E,X,0)\leq l$. Thus $l$ satisfies our requirements.
\end{proof}

\subsection{Uniform rational polytopes}

\begin{proof}[Proof of Theorem \ref{thm: uniform rational polytope foliation intro}]
The question is local, so we may work over an open neighborhood of a closed point $x\in X$. If $\pld(X\ni x,\Ff)=0$, then $B=0$ near $x$ and there is nothing left to prove. So we may assume that $\pld(X\ni x,\Ff)>0$. By Theorem \ref{thm: classification of foliation intro}, $(X\ni x,\Ff,B)$ is either smooth or as in Case 1 of Theorem \ref{thm: classification of foliation intro}. By Theorem \ref{thm: classification of foliation intro}(Case 1.d), there exists a unique prime $\Ff$-invariant curve $L$ passing through $x$, such that $(X\ni x,B+L)$ is lc. By \cite[Theorem 5.6]{HLS19}, there exists an open set $U\ni\mathbf{v}_0$ of the rational polytope of $\mathbf{v}_0$, depending only on $\mathbf{v}_0$, such that $(X\ni x,\sum_{i=1}^m v_iB_i+L)$ is lc for any $(v_1,\dots,v_m)\in U$. By Theorem \ref{thm: surface from foliation to plt pair}, $(X\ni x,\Ff,\sum_{i=1}^m v_iB_i)$ is lc for any $(v_1,\dots,v_m)\in U$. The theorem immediately follows.
\end{proof}

\section{Foliated version of some conjectures in the MMP}\label{sec: foliation conjectures}

In this section, we formally introduce the foliated version of some standard conjectures of the minimal model program and discuss their background. Since the foundations of the minimal model program for foliations in dimension $\geq 4$ has not been established, it may be too ambitious to tackle these conjectures in dimension $\geq 4$ at the moment. Nevertheless, special cases of these conjectures may still be tackable in high dimensions, e.g. algebraically intergrable foliations, or Property $(*)$ foliations (\cite[Definition 3.5]{ACSS21}).

\subsection{Complements}

\begin{conj}[Complement]\label{conj: complement conjecture foliation}
Let $\epsilon$ be a positive real number, $d$ a positive integer, and $\Ii\subset[0,1]$ a DCC set.Then there exists a positive real number $n$ depending only on $\epsilon, d$ and $\Ii$ satisfying the following.

Assume that $(X/Z\ni z,\Ff,B)$ is an $(\epsilon,\Rr)$-complementary foliated triple such that $\dim X=d$ and $B\in\Ii$. Assume that either $\epsilon=0$, or $-K_{\Ff}$ is big over $Z$. Then $(X/Z\ni z,\Ff,B)$ has an $(\epsilon,n)$-complement. Moreover, if $\bar\Ii\subset\mathbb Q$, then we $(X/Z\ni z,\Ff,B)$ has a monotonic $(\epsilon,n)$-complement.
\end{conj}

Conjecture \ref{conj: complement conjecture foliation} is an analogue of Shokurov's boundedness of $(\epsilon,n)$-complement conjecture \cite[Conjecture 6.1]{CH21}. When $\Ff=T_X$, Conjecture \ref{conj: complement conjecture foliation} is generally known when $\epsilon=0$ and $X$ is of Fano type over $Z$ (\cite{Bir19,HLS19,Sho20}) and when $\dim X=2$ \cite{CH21}. We remark that the condition ``$-K_{\Ff}$ is big over $Z$" is almost an empty condition when $Z=\{pt\}$ (cf. \cite[Theorem 5.1]{AD13}, \cite[Theorem 1.1]{Dru17}) since we have restrictions of singularities. Therefore, the interesting cases of Conjecture \ref{conj: complement conjecture foliation} should appear when either $\epsilon=0$ or $\dim Z>0$.

\smallskip

Two special cases of Conjecture \ref{conj: complement conjecture foliation} are the local index conjecture and the global index conjecture:

\begin{conj}[Local index conjecture]\label{conj: local index conjecture foliation}
Let $d$ be a positive integer, $a$ a rational number, and $\Ii\subset [0,1]\cap\mathbb Q$ a DCC set. Then there exists a positive integer $I$ depending only on $d,a$ and $\Ii$ satisfying the following. 

Assume that $(X\ni x,\Ff,B)$ is a foliated germ of dimension $d$, such that $B\in\Ii$ and $\mld(X\ni x,\Ff,B)=a$. Then $I(K_{\Ff}+B)$ is Cartier near $x$.
\end{conj}
Conjecture \ref{conj: local index conjecture foliation} is an analogue of Shokurov's local index conjecture \cite[Question 5.2]{Kaw15}. Theorem \ref{thm: surface index theorem local} proves Conjecture \ref{conj: local index conjecture foliation} when $\dim X=2$. When $\Ff=T_X$, Conjecture \ref{conj: local index conjecture foliation} is known for surfaces \cite{CH21} (by classification \cite{Sho92} when $B=0$), terminal threefolds \cite{HLL22} (by classification \cite{Kaw92} when $B=0$), canonical threefolds when $B=0$ \cite{Kaw15}, log toric pairs \cite{Amb09}, and quotient singularities when $B=0$ \cite{NS22a}.

\begin{conj}[Global index conjecture]\label{conj: global index conjecture foliation}
Let $d$ be a positive integer and $\Ii\subset [0,1]\cap\mathbb Q$ a DCC set. Then there exists a positive integer $I$ depending only on $d$ and $\Ii$ satisfying the following. 

Assume that $(X,\Ff,B)$ is a projective lc foliated triple of dimension $d$, such that $B\in\Ii$ and $K_{\Ff}+B\equiv 0$. Then $I(K_{\Ff}+B)\sim 0$.
\end{conj}
Conjecture \ref{conj: global index conjecture foliation} is an analogue of Shokurov's glocal index conjecture \cite[Conjecture 6.2]{CH21}. \cite{LLM23} proves Conjecture \ref{conj: global index conjecture foliation} when $d=3$ and $B\not=0$ or when $d=2$. When $\Ff=T_X$, Conjecture \ref{conj: local index conjecture foliation} is known for surfaces \cite{PS09} (see also \cite{Bla95,Zha91,Zha93}), threefolds \cite{Xu19} (see also \cite{Jia21}), and when $-K_X$ is big \cite{HX15} (see also \cite{Bir19}).

\begin{prop}
Conjecture \ref{conj: complement conjecture foliation} for foliations in dimension $d$ of rank $r$ implies Conjectures \ref{conj: local index conjecture foliation} and \ref{conj: global index conjecture foliation} for foliations in dimension $d$ of rank $r$.
\end{prop}
\begin{proof}
Under the setting of Conjecture \ref{conj: local index conjecture foliation}, by Conjecture \ref{conj: complement conjecture foliation}, $(X\ni x,\Ff,B)$ has a monotonic $(a,I)$-complement $(X\ni x,\Ff,B^+)$ for some $I$ depending only on $a,d$ and $\Ii$. Then $B^+=B$, so $I(K_{\Ff}+B)$ is Cartier near $x$.

Under the setting of Conjecture \ref{conj: global index conjecture foliation}, $(X,\Ff,B)$ has a monotonic $I$-complement $(X,\Ff,B^+)$ for some $I$ depending only on $d$ and $\Ii$. Thus $B^+=B$, so $I(K_{\Ff}+B)\sim 0$.
\end{proof}

It is also worth to mention the global ACC conjecture for foliations.

\begin{conj}[Global ACC]\label{conj: foliation global acc}
  Let $d$ be a positive integer and $\Ii\subset [0,1]$ a DCC set. Then there exists a finite set $\Ii_0\subset\Ii$ satisfying the following.

Assume that $(X,\Ff,B)$ is a projective lc foliated triple of dimension $d$, such that $B\in\Ii$ and $K_{\Ff}+B\equiv 0$. Then $B\in\Ii_0$.
\end{conj}
Conjecture \ref{conj: foliation global acc} is an analogue of the global ACC for usual pairs \cite[Theorem 1.5]{HMX14}.
\cite{Che22} proves Conjecture \ref{conj: foliation global acc} when $d=2$, while \cite{LLM23} proves Conjecture \ref{conj: foliation global acc} when $d=3$ and $\Ii\subset\Qq$. We remark that it is clear that Conjecture \ref{conj: global index conjecture foliation} for foliations in dimension $d$ of rank $r$ implies Conjecture \ref{conj: foliation global acc}  for foliations in dimension $d$ of rank $r$ such that $\Ii\subset\mathbb Q$.

Finally, we recall the following conjecture on the boundedness of Fano foliations:
\begin{conj}[cf. {\cite[Page 5, Problem]{Ara}}]\label{conj: bounded foliation}
Let $d$ be a positive integer. Then Fano foliations on smooth projective varieties of dimension $d$ form a bounded family. 
\end{conj}
Since the boundedness of complements \cite{Bir19} is the key to prove the BAB conjecture \cite{Bir21}, we expect \ref{conj: complement conjecture foliation} to be useful for the solution of Conjecture \ref{conj: bounded foliation}.

\subsection{Minimal log discrepancies}

We have two additional conjectures related to the mlds of foliations.

\begin{conj}[ACC for mlds]\label{conj: acc mld foliation}
       Let $d$ be a positive integer and $\Ii\subset[0,1]$ a DCC set. Then
       $$\{\mld(X\ni x,\Ff,B)\mid (X\ni x,\Ff,B) \text{ is lc}, \dim X=d,B\in\Ii\}$$
       satisfies the ACC.
\end{conj}
Conjecture \ref{conj: acc mld foliation} is an analogue of Shokurov's ACC conjecture for minimal log discrepancies \cite[Problem 5]{Sho88}. Conjecture \ref{conj: acc mld foliation} is known when $d=2$ (\cite[Theorem 0.2]{Che23}, Corollary \ref{cor: acc mld foliation surface reprove}). When $\Ff=T_X$,  Conjecture \ref{conj: acc mld foliation} is known for surfaces \cite{Ale93}, log toric pairs \cite{Amb06}, exceptional singularities \cite{HLS19}, quotient singularities when $B=0$ \cite{NS22b}, and many cases in dimension $3$ \cite{Kaw92,Mar96,Kaw15,Nak16,Jia21,LX21,HL22a,HLL22,LL22,NS22b}.

\begin{conj}[Uniform boundedness of mlds]\label{conj: bdd mld computing divisor foliation}
    Let $d$ be a positive integer and $\Ii\subset[0,1]$ a DCC set. Then there exists a positive real number $l$ depending only on $d$ and $\Ii$ satisfying the following. 
    
    Assume that $(X\ni x,\Ff,B)$ is an lc foliated germ of dimension $d$ such that $K_{\Ff}$ is $\Qq$-Cartier and $B\in\Ii$. Then there exists a prime divisor $E$ over $X\ni x$, such that $a(E,\Ff,B)=\mld(X\ni x,\Ff,B)$ and $a(E,\Ff,0)\leq l$.
\end{conj}
Conjecture \ref{conj: bdd mld computing divisor foliation} is an analogue of the uniform boundedness conjecture for mlds \cite[Conjecture 8.2]{HLL22} (see \cite[Conjecture 1.1]{MN18} for an embryonic form). Theorem \ref{thm: uniform bound of mlds} proves Conjecture \ref{conj: bdd mld computing divisor foliation} when $d=2$. When $\Ff=T_X$, Conjecture \ref{conj: bdd mld computing divisor foliation} is known for surfaces \cite{HL22b} (see \cite{MN18} for the ideal-adic case when $\Ii$ is a finite set), terminal threefold pairs \cite{HLL22}, and log toric pairs \cite{HLL22}.

Finally, we remark that we expect the generalized pair version \cite{BZ16,HL23} (more precisely, the generalized foliated quadruple version \cite{CHLX23,LLM23}) of the conjectures we have mentioned in this section to hold as well. For the reader's convenience, we omit the details of these conjectures.

\noindent\textbf{Acknowledgments}. We thank Paolo Cascini, Guodu Chen, Yen-An Chen, Christopher D. Hacon, Yuchen Liu, and Yujie Luo for helpful discussions and comments on the manuscript.  We would like to acknowledge the assistance of ChatGPT in polishing the wording. We thank A. Langer for the useful discussions on the Grauert-Riemenschneider type vanishing theorem and the classification of foliated lc surface singularities in this paper, and for providing us with the references \cite{HL21,Lup21,Lup}.

\section*{Declarations}

\begin{itemize}
\item Funding: The third author is partially supported by NSF research grants no: DMS-1801851, DMS-1952522 and by a grant from the Simons Foundation; Award Number: 256202. 
\item Conflict of interest/Competing interests: There is no conflict of interests.
\item Ethical approval: The writing of this paper is ethical. This is a pure mathematics paper so there is no experimental protocol.
\item Consent to participate: All authors consent to participate in writing this paper.
\item Consent for publication: All authors consent to participate to be published in European Journal of Mathematics, if accepted.
\item Availability of data and materials: No data.
\item Code availability: No code.
\item Authors' contributions: All authors contributed equally to this work.
\end{itemize}

\end{document}